\documentclass{jonasart}

\usepackage{enumitem}

\DeclareMathOperator{\ad}{ad}
\DeclareMathOperator{\bi}{bi}
\DeclareMathOperator{\ch}{char}
\DeclareMathOperator{\coev}{coev}
\DeclareMathOperator{\End}{End}
\DeclareMathOperator{\ev}{ev}
\DeclareMathOperator{\hemi}{hemi}
\DeclareMathOperator{\Hom}{Hom}
\DeclareMathOperator{\id}{id}

\DeclareMathOperator{\irr}{Irr}
\DeclareMathOperator{\K}{Gr}
\DeclareMathOperator{\Ker}{Ker}
\DeclareMathOperator{\leib}{Leib}
\DeclareMathOperator{\lie}{Lie}
\DeclareMathOperator{\sym}{sym}
\DeclareMathOperator{\U}{U}
\DeclareMathOperator{\UL}{UL}
\DeclareMathOperator{\weak}{weak}

\newcommand{\aaf}{\mathfrak{a}}
\newcommand{\af}{\mathfrak{A}}
\newcommand{\cbb}{\mathbb{C}}
\newcommand{\ccal}{\mathcal{C}}
\newcommand{\dcal}{\mathcal{D}}
\newcommand{\ef}{\mathfrak{e}}
\newcommand{\fbb}{\mathbb{F}}
\newcommand{\fdmod}{\mathbf{mod}}
\newcommand{\gf}{\mathfrak{g}}
\newcommand{\gl}{\mathfrak{gl}}
\newcommand{\lf}{\mathfrak{L}}
\newcommand{\Mod}{\mathbf{Mod}}
\newcommand{\nbb}{\mathbb{N}}
\newcommand{\nf}{\mathfrak{N}}
\newcommand{\ovo}{\overline{\otimes}}
\newcommand{\slf}{\mathfrak{sl}}
\newcommand{\ssf}{\mathfrak{S}}
\newcommand{\und}{\underline{\otimes}}
\newcommand{\vect}{\mathbf{Vec}}
\newcommand{\zbb}{\mathbb{Z}}

\title{Tensor products of Leibniz bimodules and Grothendieck rings}

\authors{J\"org Feldvoss and Friedrich Wagemann}

\authorinfo[J. Feldvoss]{University of South Alabama, USA}{jfeldvoss@southalabama.edu}

\authorinfo[F. Wagemann]{Nantes University, France}{wagemann@math.univ-nantes.fr}

\abstract{%
    In this paper we define three different notions of tensor products for Leibniz bimodules. The ``natural'' tensor product of Leibniz bimodules is not always a~Leibniz bimodule. In order to fix this, we introduce the notion of a~weak Leibniz bimodule and show that the ``natural" tensor product of weak bimodules is again a~weak bimodule. Moreover, it turns out that weak Leibniz bimodules are modules over a~cocommutative Hopf algebra canonically associated with the Leibniz algebra. Therefore, the category of all weak Leibniz bimodules is symmetric monoidal and the full subcategory of finite-dimensional weak Leibniz bimodules is rigid and pivotal. On the other hand, we introduce two truncated tensor products of Leibniz bimodules, which are again Leibniz bimodules. These tensor products induce a~non-associative multiplication on the Grothendieck group of the category of finite-dimensional Leibniz bimodules. In particular, we prove that in characteristic zero for a~finite-dimensional solvable Leibniz algebra over an algebraically closed field, this Grothendieck ring is an alternative power-associative commutative Jordan ring, but for a finite-dimensional non-zero semi-simple Leibniz algebra, it is neither alternative nor a~Jordan ring.
    }



\VOLUME{1}
\ISSUE{1}
\NUMBER{4}
\DOI{https://doi.org/10.46298/jonas.17482}
\editinfo{February 10, 2026}{April 14, 2026}{Salvatore Siciliano}

\acknowledgments{The authors are grateful to one of the referees for several suggestions that improved the exposition of the paper.}

\begin{document}
	\section{Introduction}

	Leibniz algebras (see~\cite{L, CH, LP1}) are non-associative algebras which generalize Lie algebras in the sense that a~Leibniz algebra has a~bilinear multiplication which is not necessarily anti-commutative, but whose left or right multiplication operators are derivations. A \textit{Leibniz
	$\lf$-bimodule\/} $M$ is a~Beck module over a~Leibniz algebra $\lf$
	described abstractly as an abelian group object in the slice category or, more concretely, as the structure we obtain on the vector space $M$
	when we declare $\lf\oplus M$ to be a~Leibniz algebra, or more explicitly, a two-sided module $M$ whose left and right $\lf$-actions are bilinear and satisfy the conditions
	\begin{gather}
		\tag{LLM}
        (xy)\cdot m=x\cdot(y\cdot m)-y\cdot(x\cdot m) \\
		\tag{LML}
        (x\cdot m)\cdot y=x\cdot(m\cdot y)-m\cdot(xy) \\
		\tag{MLL}
        (m\cdot x)\cdot y=m\cdot(xy)-x\cdot(m\cdot y)
	\end{gather}
	for all elements $x,y\in\lf$ and $m\in M$ (see~\cite{LP1} for right Leibniz algebras and~\cite{F1} for left Leibniz algebras which we will consider in this paper).

\medskip 

	The category of modules over a~Lie algebra admits a~natural tensor product over the ground field. In this paper we ask how to define a~similar tensor product for bimodules over a~Leibniz algebra $\lf$. This is not obvious as $\lf$-bimodules are (right or left) modules over the universal enveloping algebra of $\lf$ (see~\cite{LP1}), but it is not known whether the latter is a~bialgebra, and our paper seems to indicate that this might not be the case. It turns out that the tensor product $M\otimes N$ of $\lf$-bimodules $M$ and $N$ endowed with the left $\lf$-action
	\[
		x\cdot(m\otimes n)=(x\cdot m)\otimes n+m\otimes(x\cdot n)
	\]
	and the right $\lf$-action
	\[
		(m\otimes n)\cdot x=(m\cdot x)\otimes n+m\otimes(n\cdot x)
	\]
	satisfies the conditions (LLM) and (LML), but in general not the condition (MLL) (Proposition~\ref{tensprod} and Example~\ref{weak1dim}).
\medskip 

	Based on this result, we pursue two lines of thoughts to define a~tensor product on some category of bimodules over a~Leibniz algebra. First of all, in Section~\ref{weak}, we introduce the concept of a~\textit{weak Leibniz bimodule\/} as a two-sided module over a~Leibniz algebra satisfying (LLM) and (LML), but not necessarily (MLL). Then we show that weak Leibniz bimodules are the left modules over a~cocommutative Hopf algebra canonically associated with the Leibniz algebra, its \textit{weak universal enveloping algebra\/} (Theorem~\ref{hopf}). We deduce from here that the category $\Mod^{\bi}_{\weak}(\lf)$ of weak $\lf$-bimodules is a~symmetric monoidal category (Theorem~\ref{monoidal}) and that its full subcategory $\fdmod^{\bi}_{\weak}(\lf)$ of finite-dimensional weak $\lf$-bimodules is rigid and pivotal (Theorem~\ref{rigid} and Remark~\ref{pivotal}). In fact, $\fdmod^{\bi}_{\weak}(\lf)$ is even a~ring category (Theorem~\ref{ring}) in the sense of~\cite{EGNO}. Along the way, we address the question of classifying the irreducible objects in $\Mod^{\bi}_{\weak}(\lf)$ and obtain a~partial answer (Proposition~\ref{weakirred}), we discuss the relation between the category $\Mod^{\bi}_{\weak}(\lf)$ and the full subcategory $\Mod^{\bi}(\lf)$ of all $\lf$-bimodules (Proposition~\ref{reflectiv}), and we give some examples designed to illustrate our results.
\medskip 

	The second line of thought is to consider the tensor product $M\otimes N$ of the Leibniz bimodules $M$ and $N$ modulo a~certain subspace in order to obtain an ordinary Leibniz bimodule. We look at two such subspaces which thus lead to two \textit{truncated\/} tensor products $M\ovo N$ and $M\und N$ whose properties are discussed in Section~\ref{truncated}. These truncated tensor products are the same when one tensor factor is symmetric or anti-symmetric (Theorem~\ref{main}), and they coincide with the ``natural'' tensor product defined in Section~\ref{weak} in the case that both factors are symmetric or anti-symmetric (Corollary~\ref{sym}). On the other hand, the truncated tensor products are zero when both factors are non-trivial irreducible, one factor is symmetric and the other one is anti-symmetric (Corollary~\ref{irred}). We give examples to show that in many cases both truncated tensor products are non-zero. But an important property which one looses by considering the truncated tensor products is the associativity.
\medskip 

	In order to capture this loss of associativity more precisely, we construct in Section~\ref{grothendieck} the Grothendieck ring associated with the truncated tensor products and describe its algebraic structure. Note that we do not know whether $\ovo$ and $\und$ do coincide, but they always induce the same multiplication on the Grothendieck group $\K^{\bi}(\lf)$ of the category $\fdmod^{\bi}(\lf)$ of finite-dimensional $\lf$-bimodules for a~Leibniz algebra $\lf$. This Grothendieck ring is a~unital commutative ring (Proposition~\ref{unitalcomm}), but it is not necessarily associative. One of the main results in this section, Theorem~\ref{communitalprod}, is the description of the Grothendieck ring $\K^{\bi}(\lf)$ as a~so-called \textit{unital commutative product\/} of two copies of the Grothendieck ring $\K(\lf_{\lie})$ of the canonical Lie algebra $\lf_{\lie}$ associated with the Leibniz algebra $\lf$. Here one copy of $\K(\lf_{\lie})$ corresponds to the classes of the symmetric irreducible $\lf$-bimodules and the other copy corresponds to the classes of the anti-symmetric irreducible $\lf$-bimodules each of which can be constructed in a~natural way from the irreducible $\lf_{\lie}$-modules. We use this result to completely determine the structure of $\K^{\bi}(\lf)$ for a~finite-dimensional solvable Leibniz algebra $\lf$ over an algebraically closed field of characteristic zero (Corollary~\ref{solv}), and we also prove that in this case $\K^{\bi}(\lf)$ is an alternative power-associative Jordan ring (Theorem~\ref{alt} and Corollary~\ref{alt2}). In addition, we give an example to show that Theorem~\ref{alt} is not true for ground fields of non-zero characteristic (Example~\ref{primchar}). On the other hand, we prove in Proposition~\ref{nonassoK_0} that $\K^{\bi}(\lf)$ is not associative for those Leibniz algebras $\lf$ whose canonical Lie algebras $\lf_{\lie}$ are finite dimensional. More specifically, we show that the Grothendieck ring $\K^{\bi}(\lf)$ of a~finite-dimensional non-zero semi-simple Leibniz algebra $\lf$ over a~field of characteristic zero is neither alternative nor a~Jordan ring (Theorem~\ref{nonalt}).
\medskip 

	Although we are able to obtain several results for (truncated) tensor products of
	Leibniz bimodules, many questions still remain open. Contrary to irreducible Leibniz bimodules, a~classification of the irreducible weak Leibniz bimodules (up to isomorphism)
	is not known. Such a~classification might also shed some light on the precise relationship between the categories $\Mod^{\bi}_{\weak}(\lf)$ and $\Mod^{\bi}(\lf)$ of (weak)
	$\lf$-bimodules or the Grothendieck rings of the corresponding subcategories
	$\fdmod^{\bi}_{\weak}(\lf)$ and $\fdmod^{\bi}(\lf)$ of finite-dimensional (weak)
	$\lf$-bimodules. Moreover, except for finite-dimensional solvable Leibniz algebras over an algebraically closed field of characteristic zero and $\slf_2 (\cbb)$, we do not know whether the Grothendieck ring of the category of finite-dimensional
	Leibniz bimodules is power-associative. Finally, it would be very interesting to explicitly compute the Grothendieck rings for non-solvable Leibniz algebras whose canonical Lie algebra is not isomorphic to $\slf_2 (\cbb)$ or for Leibniz algebras over a~field of non-zero characteristic.
\medskip 

	We finish this introduction by fixing some notation and conventions. A ring without any specification is an abelian group with a~biadditive multiplication that not necessarily satisfies any other identities. Similarly, an algebra without any specification is a~vector space with a~bilinear multiplication that not necessarily satisfies any other identities. The multiplicative identity element of a~ring or an algebra is said to be its unity, a~ring or an algebra with a~unity is called unital, and homomorphisms between unital rings or unital algebras preserve the unity. Ideals are always two-sided ideals if not explicitly stated otherwise. Modules over Lie algebras are always left modules. All vector spaces and algebras are defined over an arbitrary field which is only explicitly mentioned when some additional assumptions on the ground field are made or this enhances the understanding of the reader. In particular, with one exception, when we consider tensor products over the integers, which will be denoted by $\otimes_\zbb$, all tensor products are over the relevant ground field and will be denoted by $\otimes$. For a~subset $X$ of a~vector space $V$ over a~field $\fbb$ we let $\langle X\rangle_\fbb$ be the subspace of $V$ spanned by $X$. We denote the space of linear transformations from an $\fbb$-vector space $V$ to an $\fbb$-vector space $W$ by $\Hom_\fbb(V,W)$. In particular, $\End_\fbb(V):=\Hom_\fbb(V,V)$ is the space of linear operators on $V$, and the linear dual $V^*:=\Hom_\fbb(V,\fbb)$
	is the space of linear forms on $V$. Finally, let $\id_X$ denote the identity function on a~set $X$, let $\nbb_0$ be the set of non-negative integers, and let $\zbb$ be the ring of integers.
\newpage  

	\section{Preliminaries}

	In this section we recall some definitions and prove a~result about the universal enveloping algebra of a~left Leibniz algebra that will be useful later in the paper.

	A \emph{left Leibniz algebra\/} is an algebra $\lf$ such that every left multiplication operator $L_x:\lf\to\lf$, $y\mapsto xy$ is a~derivation. This is equivalent to the identity
	\[
		x(yz)=(xy)z+y(xz)
	\]
	for all elements $x,y,z\in\lf$, which in turn is equivalent to the identity
	\[
		(xy)z=x(yz)-y(xz)
	\]
	for all elements $x,y,z\in\lf$. There is a~similar definition of a~\emph{right Leibniz algebra\/}, which is used by Loday et al.\ (see~\cite{L, CH, LP1}), but as in our previous papers (see~\cite{F1, F2, FW1}), we will consider left Leibniz algebras which most of the time are just called Leibniz algebras.

	Note that every Lie algebra is a~left and a~right Leibniz algebra. On the other hand, every Leibniz algebra has an important ideal, its Leibniz kernel, that measures how much the Leibniz algebra deviates from being a~Lie algebra. Namely, let $\lf$ be a
	Leibniz algebra over a~field $\fbb$. Then
	\[
		\leib(\lf):=\langle x^2\mid x\in\lf\rangle_\fbb
	\]
	is called the \emph{Leibniz kernel\/} of $\lf$. The Leibniz kernel $\leib(\lf)$ is an abelian ideal of $\lf$, and $\leib(\lf)\ne\lf$ whenever $\lf\ne 0$ (see~\cite[Proposition 2.20]{F1}). Moreover, $\lf$ is a~Lie algebra if, and only if, $\leib(\lf)=0$.

	By definition of the Leibniz kernel, $\lf_{\lie}:=\lf/\leib(\lf)$ is a~Lie algebra which we call the \emph{canonical Lie algebra\/} associated with $\lf$. In fact, the Leibniz kernel is the smallest ideal such that the corresponding factor algebra is a~Lie algebra (see~\cite[Proposition~2.22]{F1}).

	Next, we will briefly discuss left modules and bimodules of left Leibniz algebras. Let
	$\lf$ be a~left Leibniz algebra over a~field $\fbb$. A \emph{left $\lf$-module\/} is a vector space $M$ over $\fbb$ with an $\fbb$-bilinear left $\lf$-action $\lf\times M\to
	M$, $(x,m)\mapsto x\cdot m$ such that
	\[
		(xy)\cdot m=x\cdot(y\cdot m)-y\cdot(x\cdot m)
	\]
	is satisfied for every $m\in M$ and all $x,y\in\lf$. Note that for a~Lie algebra a~left
	Leibniz module is the same as a~Lie module.

	Moreover, every left $\lf$-module $M$ gives rise to a~homomorphism $\lambda:\lf\to
	\gl(M)$ of left Leibniz algebras, defined by $\lambda_x (m):=x\cdot m$, and vice versa. We call $\lambda$ the \emph{left representation\/} of $\lf$ associated with $M$.

	The correct concept of a~module for a~left Leibniz algebra $\lf$ is the notion of a~Leibniz
	$\lf$-bimodule. An \emph{$\lf$-bimodule\/} is a~vector space $M$ that is endowed with an $\fbb$-bilinear left $\lf$-action and an $\fbb$-bilinear right $\lf$-action satisfying the following compatibility conditions:
		\begin{gather}
		\tag{LLM}
        (xy)\cdot m=x\cdot(y\cdot m)-y\cdot(x\cdot m) \\
		\tag{LML}
        (x\cdot m)\cdot y=x\cdot(m\cdot y)-m\cdot(xy) \\
		\tag{MLL}
        (m\cdot x)\cdot y=m\cdot(xy)-x\cdot(m\cdot y)
	\end{gather}
    for all elements $x,y\in\lf$ and $m\in M$.

	It is an immediate consequence of (LLM) that every Leibniz bimodule is a~left Leibniz module.
	\begin{remark}\label{ZD}
		{Note that when (LML) holds, the relation (MLL) is equivalent to
		\begin{equation}
			\tag{ZD}
            (x\cdot m+m\cdot x)\cdot y=0
		\end{equation}
		for all elements $x,y\in\lf$ and $m\in M$.}
	\end{remark}
	On the other hand, a~pair $(\lambda,\rho)$ of linear transformations $\lambda:\lf\to
	\End_\fbb(V)$ and $\rho:\lf\to\End_\fbb(V)$ is called a~\emph{representation\/} of $\lf$
	on the $\fbb$-vector space $V$ if the following conditions are satisfied:
	\begin{gather}
		\lambda_{xy}=\lambda_x\circ\lambda_y -\lambda_y\circ\lambda_x \label{LLMrep} \\
		\rho_{xy}=\lambda_x\circ\rho_y -\rho_y\circ\lambda_x \label{LMLrep} \\
		\rho_y\circ(\lambda_x +\rho_x )=0 \label{MLLrep}
	\end{gather}
	for all elements $x,y\in\lf$. Then every $\lf$-bimodule $M$ gives rise to a~representation
	$(\lambda,\rho)$ of $\lf$ on $M$ via $\lambda_x (m):=x\cdot m$ and $\rho_x (m):=m\cdot x$. Conversely, every representation $(\lambda,\rho)$ of $\lf$ on the vector space $M$
	defines an $\lf$-bimodule structure on $M$ via $x\cdot m:=\lambda_x (m)$ and $m\cdot x:=\rho_x (m)$.

	By virtue of~\cite[Lemma 3.3]{F1}, every left $\lf$-module is an $\lf_{\lie}$-module in a natural way, and vice versa. Consequently, many properties of left Leibniz modules follow from the corresponding properties of modules for the canonical Lie algebra.

	The usual definitions of the notions of \emph{sub\textrm{(}bi\textrm{)}module\/},
	\emph{irreducibility\/}, \emph{complete reducibility\/}, \emph{composition series\/},
	\emph{homomorphism\/}, \emph{isomorphism\/}, etc., hold for left Leibniz modules and Leibniz bimodules. (Note that by definition an irreducible (bi)module is always non-zero.)

	Let $M$ be an $\lf$-bimodule for a~left Leibniz algebra $\lf$. Then $M$ is said to be \emph{symmetric\/} if $m\cdot x=-x\cdot m$ for every $x\in\lf$ and every
	$m\in M$, and $M$ is said to be \emph{anti-symmetric\/} if $m\cdot x=0$ for every $x\in\lf$ and every $m\in M$. Moreover, an $\lf$-bimodule $M$ is called
	\emph{trivial\/} if $x\cdot m=0=m\cdot x$ for every $x\in\lf$ and every $m\in M$. Note that an $\lf$-bimodule $M$ is trivial if, and only if, $M$ is symmetric and anti-symmetric. We call
	\[
		M_0:=\langle x\cdot m+m\cdot x\mid x\in\lf,m\in M\rangle_\fbb
	\]
	the \emph{anti-symmetric kernel\/} of $M$. It is well known that $M_0$ is an anti-symmetric $\lf$-subbimodule of $M$ such that $M_{\sym}:=M/M_0$ is symmetric (see~\cite[Propositions~3.12 and~3.13]{F1}). Recall that every left $\lf$-module $M$ of a~left Leibniz algebra $\lf$ determines a unique symmetric $\lf$-bimodule structure on $M$ by defining $m\cdot x:=-x\cdot m$ for every element $m\in M$ and every element $x\in\lf$
	(see~\cite[Proposition~3.15\,(b)]{F1}). We call this $\lf$-bimodule the
	\emph{symmetrization\/} of $M$ which will be denoted by $M^s$. Similarly, every left $\lf$-module $M$ with trivial right action is an anti-symmetric
	$\lf$-bimodule (see~\cite[Proposition~3.15\,(a)]{F1}). We call this
	$\lf$-bimodule the \emph{anti-symmetrization\/} of $M$ which will be denoted by $M^a$.

	Loday and Pirashvili~\cite[Section 2]{LP1} define a~universal enveloping algebra of a~right Leibniz algebra. We will briefly recall the analogous construction for left Leibniz algebras. For any left Leibniz algebra $\lf$
	over a~field $\fbb$ we denote by $\lf^\ell$ an isomorphic copy of the underlying $\fbb$-vector space of $\lf$ (by the isomorphism $\ell:\lf\to
	\lf^\ell$), and we denote by $\lf^r$ another isomorphic copy of the underlying
	$\fbb$-vector space of $\lf$ (by the isomorphism $r:\lf\to\lf^r$). Then the tensor algebra $T(\lf^\ell\oplus\lf^r )$ is a~unital associative algebra over $\fbb$. We would like the multiplication on the factor algebra
	$T(\lf^\ell\oplus\lf^r )/I(\lf)$ to satisfy the following relations:
	\begin{gather}
		\ell_{xy}=\ell_x\cdot\ell_y -\ell_y\cdot\ell_x \label{llm} \\
		r_{xy}=\ell_x\cdot r_y -r_y\cdot\ell_x \label{lml} \\
		r_x\cdot(\ell_y +r_y )=0 \label{zd}
	\end{gather}
	for all elements $x,y\in\lf$, where $\cdot$ denotes the usual multiplication in the tensor algebra $T(\lf^\ell\oplus\lf^r )$, or equivalently, the two-sided ideal $I(\lf)$ of $T(\lf^\ell\oplus\lf^r )$ is defined as the two-sided ideal generated by
	\[
		\{\ell_x\cdot\ell_y -\ell_y\cdot\ell_x -\ell_{xy}\,,\,\ell_x\cdot r_y
		-r_y\cdot\ell_x -r_{xy}\,,r_x\cdot(\ell_y +r_y )\mid x,y\in\lf\}\,.
	\]
	Then $\UL(\lf):=T(\lf^\ell\oplus\lf^r )/I(\lf)$ is called the \emph{universal enveloping algebra\/} of $\lf$.
	\begin{remark}\label{rel}
		{Note that the relations \eqref{llm}, \eqref{lml}, and \eqref{zd}
		correspond to the identities \eqref{LLMrep}, \eqref{LMLrep}, and
		\eqref{MLLrep}, respectively.}
	\end{remark}
	It follows from Remark~\ref{rel} that every $\lf$-bimodule is a~unital left $\UL(\lf)$-module and vice versa (see~\cite[Theorem 2.3]{LP1}
	for the analogue for right Leibniz algebras). In other words, we have an equivalence of categories:
	\[
		\Mod^{\bi}(\lf)\cong\UL(\lf)-\Mod\,,
	\]
	where $\Mod^{\bi}(\lf)$ denotes the category of $\lf$-bimodules and
	$\UL(\lf)-\Mod$ denotes the category of unital left $\UL
	(\lf)$-modules.

	For right Leibniz algebras the following result is due to Loday and
	Pirashvili~\cite[Proposition 2.5]{LP1}:
	\begin{proposition}\label{augmented}
		Let $\lf$ be a~left Leibniz algebra over a~field $\fbb$. Then the following statements hold:
		\begin{enumerate}
			\item[\textrm{(a)}] The function $d_0:\UL(\lf)\to\U(\lf_{\lie})$ defined by
			\[
                d_0 (1):=1, \quad
                d_0 (\ell_x ):=\overline{x} \quad \text{and} \quad
                d_0 (r_x ):=0
            \]
            is a~homomorphism of unital associative $\fbb$-algebras.

			\item[\textrm{(b)}] The function $d_1:\UL(\lf)\to\U(\lf_{\lie})$ defined by
			\[
                d_1 (1):=1, \quad
                d_1 (\ell_x ):=\overline{x} \quad \text{and} \quad
                d_1 (r_x ):=-\overline{x}
            \]
            is a~homomorphism of unital associative $\fbb$-algebras.

			\item[\textrm{(c)}] The function $s_0:\U(\lf_{\lie})\to\UL(\lf)$ defined by
			$s_0 (1):=1$ and $s_0 (\overline{x}):=\ell_x$ is a homomorphism of unital associative $\fbb$-algebras.
			\item[\textrm{(d)}] $d_0\circ s_0 =\id_{\U(\lf_{\lie})}$ and $d_1\circ s_0
			=\id_{\U(\lf_{\lie})}$. In particular, $d_0$ and $d_1$
			are surjective as well as $s_0$ is injective.
			\item[\textrm{(e)}] $\Ker(d_0 )\Ker(d_1 )=0$.
		\end{enumerate}
	\end{proposition}

	\begin{proof}
		In this proof we let $I_0 (\lf_{\lie})$ denote the ideal in the tensor algebra $T(\lf_{\lie})$ generated by $\{x\otimes y-y\otimes x-xy\mid x,y\in\lf_{\lie}\}$ and set $\overline{t}:=t+I_0 (\lf_{\lie})$ for every element $t\in T(\lf_{\lie})$. In particular, we have that $\U(\lf_{\lie})
		=T(\lf_{\lie})/I_0 (\lf_{\lie})$ whose multiplication is denoted by $\odot$.

		(a): Since $d_0$ is defined on generators of the algebra $\UL(\lf)$, we only need to show that $d_0$ is well-defined, i.e., $d_0 [I(\lf)]
		=0$. Namely, we have that
		\begin{eqnarray*}
			d_0 (\ell_x\cdot\ell_y -\ell_y\cdot\ell_x -\ell_{xy}) & = &
			d_0 (\ell_x )\odot d_0 (\ell_y )-d_0 (\ell_y )\odot d_0 (\ell_x )-d_0 (\ell_{xy})\\
			& = & \overline{x}\odot\overline{y}-\overline{y}\odot
			\overline{x}-\overline{xy}\\
			& = & \overline{x\otimes y-y\otimes x-xy}=0\,,
		\end{eqnarray*}

		\begin{eqnarray*}
			d_0 (\ell_x\cdot r_y -r_y\cdot\ell_x -r_{xy}) & = & d_0 (\ell_x )
			\odot d_0 (r_y )-d_0 (r_y )
			\odot d_0 (\ell_x )-d_0 (r_{xy})\\
			& = & \overline{x}\odot 0-0\odot\overline{x}-0=0\,,
		\end{eqnarray*}
        and
		\begin{eqnarray*}
			d_0 (r_x\cdot(\ell_y +r_y )) & = & d_0 (r_x )\odot[d_0 (\ell_y )+d_0 (r_y )]\\
			& = & 0\odot(\overline{y}+0)=0\,.
		\end{eqnarray*}
		(b): As in part (a), we only need to show that $d_1 [I(\lf)]=0$. For the relation \eqref{llm} the proof is exactly the same as in part (a). For the other two relations we obtain that
		\begin{eqnarray*}
			d_1 (\ell_x\cdot r_y -r_y\cdot\ell_x -r_{xy}) & = & d_1 (\ell_x )\odot d_1 (r_y )-d_1 (r_y )\odot d_1 (\ell_x )-d_1 (r_{xy})\\
			& = & \overline{x}\odot(-\overline{y})-(-\overline{y})\odot\overline{x}
			-(-\overline{xy})\\
			& = & -(\overline{x}\odot\overline{y})+\overline{y}\odot\overline{x}
			+\overline{xy}\\
			& = & -\,\,\overline{x\otimes y-y\otimes x-xy}=0\,,
		\end{eqnarray*}
		and
		\begin{eqnarray*}
			d_1 (r_x\cdot(\ell_y +r_y )) & = & d_1 (r_x )\odot[d_1 (\ell_y )+d_1 (r_y )]\\
			& = & -\,\,\overline{x}\odot[\overline{y}+(-\overline{y})]=0\,.
		\end{eqnarray*}
		(c): Again, we have to show that $s_0$ is well-defined, i.e., that $\ell_x =0$
		for every element $x\in\leib(\lf)$, or equivalently, $\ell_{x^2}=0$ in
		$\UL(\lf)$ for every element $x\in\lf$. But the latter is an immediate consequence of relation \eqref{llm} in $\UL(\lf)$ as
		\[
			\ell_{x^2}=\ell_x\cdot\ell_x -\ell_x\cdot\ell_x =0
		\]
		for every element $x\in\lf$.

		(d): For every element $x\in\lf$ we have that
		\[
			(d_0\circ s_0 )(1)=d_0 (s_0 (1))=d_0 (1)=1
		\]
		as well as
		\[
			(d_0\circ s_0 )(\overline{x})=d_0 (s_0 (\overline{x}))=d_0 (\ell_x )=
			\overline{x}\,,
		\]
		and a~similar argument proves the second identity.

		(e): As an ideal, $\Ker(d_0 )$ is generated by $\{r_x\mid x\in\lf\}$
		and $\Ker(d_1 )$ is generated by $\{\ell_y +r_y\mid y\in\lf\}$. Hence, the assertion is an immediate consequence of the relation \eqref{zd}
		in $\UL(\lf)$.
	\end{proof}

	It is well known that the universal enveloping algebra of a~Lie algebra is a~Hopf algebra, but is not known whether the universal enveloping algebra of a~Leibniz algebra is a~Hopf algebra, and this paper seems to indicate that this might not be the case. On the other hand, it follows from Proposition~\ref{augmented} that
    \[
        \varepsilon_{\lie}\circ d_0
        = \varepsilon_{\lie}\circ d_1
        \colon \UL(\lf)\to\fbb
    \]
    is an epimorphism of unital associative $\fbb$-algebras, where $\varepsilon_{\lie}$ denotes the counit of $\U(\lf_{\lie})$. Hence, $\UL(\lf)$ is an \emph{augmented unital associative $\fbb$-algebra\/} with \emph{augmentation map\/}
	$\varepsilon_{\lie}\circ d_0 =\varepsilon_{\lie}\circ d_1$.

	\section{Weak Leibniz bimodules}\label{weak}

	We would like to make the tensor product of Leibniz bimodules into a Leibniz bimodule. Let $\lf$ be a~left Leibniz algebra over a~field
	$\fbb$, and let $M$ and $N$ be left $\lf$-modules. Then, as for modules over a~Lie algebra, the tensor product $M\otimes N$ (over the ground field $\fbb$) is a~left $\lf$-module endowed with a~left action
	\[
		x\cdot(m\otimes n):=(x\cdot m)\otimes n+m\otimes(x\cdot n)
	\]
	for all elements $x\in\lf$, $m\in M$, $n\in N$.

	This motivates the following definition of a~left and right action on the tensor product of Leibniz bimodules. For $\lf$-bimodules
	$M$ and $N$ we define the following left and right actions of a~left
	Leibniz algebra $\lf$ on the vector space $M\otimes N$:
	\[
		x\cdot(m\otimes n):=(x\cdot m)\otimes n+m\otimes(x\cdot n)
	\]
	and
	\[
		(m\otimes n)\cdot x:=(m\cdot x)\otimes n+m\otimes(n\cdot x)
	\]
	for all elements $x\in\lf$, $m\in M$, $n\in N$. Then, maybe a bit surprisingly, we have the following result:
	\begin{proposition}\label{tensprod}
		Let $\lf$ be a~left Leibniz algebra, and let $M$ and $N$ be $\lf$-bimodules. Then $M\otimes N$ satisfies \textup{(LLM)} and \textup{(LML)}. Furthermore,
		\textup{(MLL)} holds if, and only if,
		\begin{equation}\label{*}
			(x\cdot m+m\cdot x)\otimes (n\cdot y)+(m\cdot y)\otimes
			(x\cdot n+n\cdot x)=0
		\end{equation}
		holds for all elements $x,y\in\lf$, $m\in M$, and $n\in N$. In particular, if $M$ and $N$ are both symmetric or both anti-symmetric, then $M\otimes N$ is a~symmetric or anti-symmetric $\lf$-bimodule, respectively.
	\end{proposition}

	\begin{proof}
		As for the tensor product of Lie modules, $M\otimes N$ satisfies (LLM). Next, let us verify that (LML) holds for $M\otimes N$. Namely, as (LML)
		is satisfied for $M$ and $N$, we have that
		\begin{eqnarray*}
			(m\otimes n)\cdot(xy) & = & [m\cdot(xy)\otimes n]+m\otimes[n\cdot(xy)]\\
			& = & [x\cdot(m\cdot y)]\otimes n-[(x\cdot m)\cdot y]\otimes n\\
			&& +m\otimes[x\cdot(n\cdot y)]-m\otimes[(x\cdot n)\cdot y]\\
			& = & [x\cdot(m\cdot y)]\otimes n+(m\cdot y)\otimes(x\cdot n)\\
			&& +(x\cdot m)\otimes(n\cdot y)+m\otimes[x\cdot(n\cdot y)]\\
			&& -[(x\cdot m)\cdot y]\otimes n-(x\cdot m)\otimes(n\cdot y)\\
			&& -(m\cdot y)\otimes(x\cdot n)-m\otimes[(x\cdot n)\cdot y]\\
			& = & x\cdot[(m\cdot y)\otimes n+m\otimes(n\cdot y)]\\
			&& -[(x\cdot m)\otimes n+m\otimes(x\cdot n)]\cdot y\\
			& = & x\cdot[(m\otimes n)\cdot y]-[x\cdot(m\otimes n)]\cdot y
		\end{eqnarray*}
		for all elements $x,y\in\lf$, $m\in M$, and $n\in N$.

		Finally, instead of (MLL), we verify (ZD) for the tensor product
		$M\otimes N$ (see Remark~\ref{ZD}). Because (ZD) holds for $M$ and $N$, we obtain that
		\begin{align*}
			[x\cdot(m\otimes n)+(m\otimes n)\cdot x]
			\cdot y
            ={ }&{ } [(x\cdot m)\otimes n)
			+m\otimes(x\cdot n)\\
			& +(m\cdot x)\otimes n+m\otimes(n\cdot x)]\cdot y\\
			={ }&{ }[(x\cdot m)\cdot y]\otimes n+(x\cdot m)\otimes(n\cdot y)\\
			& +(m\cdot y)\otimes(x\cdot n)+m\otimes[(x\cdot n)\cdot y]\\
			& +[(m\cdot x)\cdot y]\otimes n+(m\cdot x)\otimes(n\cdot y)\\
			& +(m\cdot y)\otimes(n\cdot x)+m\otimes[(n\cdot x)\cdot y]\\
			={ }&{ }[\underbrace{(x\cdot m+m\cdot x)\cdot y}_{=0}]\otimes n
			+(x\cdot m+m\cdot x)\otimes(n\cdot y)\\
			& +(m\cdot y)\otimes(x\cdot n+n\cdot x)+m\otimes
			[\underbrace{(x\cdot n+n\cdot x)\cdot y}_{=0}]\\
			={ }&{ }(x\cdot m+m\cdot x)\otimes(n\cdot y)+(m\cdot y)\otimes(x\cdot n+n\cdot x)
		\end{align*}
		for all elements $x,y\in\lf$, $m\in M$, and $n\in N$. This shows that (ZD), or equivalently (MLL), holds exactly when
		\[
			(x\cdot m+m\cdot x)\otimes(n\cdot y)+(m\cdot y)\otimes(x\cdot n+n\cdot x)=0
		\]
		for all elements $x,y\in\lf$, $m\in M$, and $n\in N$, as asserted. Moreover, the last statement can be read off from this identity.
	\end{proof}

	\begin{remark}\label{rel2}
		{It should be noted that in the above proof of (LLM) for the tensor product of Leibniz bimodules, we only used that (LLM) is satisfied for each of the factors. Similarly, the proof of (LML) for the tensor product of Leibniz bimodules only uses that (LML) is satisfied for each of the factors.}
	\end{remark}
	In the light of Proposition~\ref{tensprod} and Remark~\ref{rel2}, we say that a left Leibniz module with a~right action satisfying (LML) is a~\emph{weak Leibniz bimodule\/}. Then the proof of Proposition~\ref{tensprod} shows that the tensor product of weak Leibniz modules is again a~weak Leibniz bimodule:
	\begin{proposition}
		Let $\lf$ be a~Leibniz algebra. If $M$ and $N$ are weak $\lf$-bimodules, then $M\otimes N$ is a~weak $\lf$-bimodule.
	\end{proposition}
	If we use the same definition of (anti-)symmetry for weak Leibniz bimodules as for ordinary Leibniz bimodules, then the following result is an immediate consequence of identity (ZD) in Remark~\ref{ZD}:
	\begin{proposition}\label{weakantisymm}
		Let $\lf$ be a~Leibniz algebra. Then the following statements hold:
		\begin{enumerate}
			\item[\textrm{(a)}] Every symmetric weak $\lf$-bimodule is an $\lf$-bimodule.
			\item[\textrm{(b)}] Every anti-symmetric weak $\lf$-bimodule is an $\lf$-bimodule.
		\end{enumerate}
	\end{proposition}
	The next result is an immediate consequence of~\cite[Theorem 3.14]{F1}
	and Proposition~\ref{weakantisymm}:
	\begin{corollary}\label{weakantisymmirred}
		Let $\lf$ be a~Leibniz algebra. Then an irreducible weak $\lf$-bimodule is an $\lf$-bimodule if, and only if, it is symmetric or anti-symmetric.
	\end{corollary}

	\begin{remark}
	   As in Proposition~\ref{tensprod}, we obtain from Proposition~\ref{weakantisymm} that the tensor product of two symmetric or two anti-symmetric weak Leibniz bimodules is a~symmetric or an anti-symmetric Leibniz bimodule, respectively. On the other hand, the tensor product of a~symmetric Leibniz bimodule and an anti-symmetric Leibniz bimodule (in any order) is in general only a weak Leibniz bimodule (see Example~\ref{weak1dim} below). In particular, this shows that the tensor product of two ordinary
		Leibniz bimodules is not always an ordinary Leibniz bimodule.
	\end{remark}
	Let $J(\lf)$ denote the two-sided ideal of the tensor algebra $T(\lf^\ell
	\oplus\lf^r )$ generated by
	\[
		\{\ell_x\cdot\ell_y -\ell_y\cdot\ell_x -\ell_{xy}\,,\,\ell_x\cdot r_y -r_y\cdot
		\ell_x -r_{xy}\mid x,y\in\lf\}\,.
	\]
	Then we say that $\UL_{\weak}(\lf):=T(\lf^\ell\oplus\lf^r )/J(\lf)$ is the
	\emph{weak universal enveloping algebra\/} of $\lf$.

	Let $\Mod_{\weak}^{\bi}(\lf)$ denote the category of weak $\lf$-bimodules for a~Leibniz algebra $\lf$, and let $\UL_{\weak}(\lf)-\Mod$
	denote the category of unital left $\UL_{\weak}(\lf)$-modules. It follows from Remark~\ref{rel} that every weak $\lf$-bimodule is a~unital left
	$\UL_{\weak}(\lf)$-module and vice versa, and thus we have an equivalence of categories:
	\[
		\Mod_{\weak}^{\bi}(\lf)\cong\UL_{\weak}(\lf)-\Mod\,.
	\]
	As a~consequence, we obtain that the category of weak Leibniz bimodules is an $\fbb$-linear abelian category (see~\cite[Chapter
	VIII]{McL} and~\cite[Definition 1.2.2]{EGNO}) as was already observed in~\cite[Section 5]{L} for the category of bimodules over a~right Leibniz algebra:
	\begin{proposition}\label{ab}
		Let $\lf$ be a~Leibniz algebra over a~field $\fbb$. Then
		$\Mod_{\weak}^{\bi}(\lf)$ is an $\fbb$-linear abelian category.
	\end{proposition}
	Next, we briefly discuss the irreducible objects of $\Mod_{\weak}^{\bi}
	(\lf)$. Recall that every irreducible $\lf$-bimodule is either symmetric or anti-sym\-metric (see~\cite[Theorem~3.14]{F1}). The main ingredient in the proof of this result is the anti-symmetric kernel $M_0$ of an
	$\lf$-bimodule $M$ which is an $\lf$-subbimodule. Consequently, either $M_0 =0$, i.e., $M$ is symmetric, or $M=M_0$, i.e., $M$ is anti-symmetric. In particular, every irreducible $\lf$-bimodule arises from an irreducible $\lf_{\lie}$-module by symmetrization or anti-symmetrization
	(see~\cite[Proposition 3.15]{F1}), and therefore the classification of the irreducible $\lf$-bimodules reduces to the classification of the irreducible modules over its canonical Lie algebra $\lf_{\lie}$. For weak
	Leibniz bimodules the situation is more complicated as the analogue of the anti-symmetric kernel is not necessarily anti-symmetric, i.e., its right action is not necessarily trivial (see Example~\ref{weak1dim}
	below). Of course, the reason for this is that weak Leibniz bimodules do not necessarily satisfy the identity (ZD). In fact, in general it is not clear to us whether the anti-symmetric kernel of a~weak Leibniz bimodule is invariant under the right action, but the proof of~\cite[Proposition~3.12]{F1}
	shows that it is still invariant under the left action.

	As a~replacement for the anti-symmetric kernel of a~weak $\lf$-bimodule
	$M$ we consider
	\[
		M\lf:=\langle m\cdot x\mid m\in M\,,\,x\in\lf\rangle_\fbb\,,
	\]
	which indeed is an $\lf$-subbimodule of $M$. Namely, it follows from
	(LML) that $M\lf$ is invariant under the left $\lf$-action and by definition
	$M\lf$ is invariant under the right $\lf$-action. Similarly, the space of right
	$\lf$-invariants
	\[
		M^\lf:=\{m\in M\mid\forall\,x\in\lf\,:\,m\cdot x=0\}
	\]
	is an $\lf$-subbimodule of $M$. Namely, it again follows from (LML)
	that $M^\lf$ is invariant under the left $\lf$-action and by definition
	$M^\lf$ is invariant under the right $\lf$-action.
	\begin{lemma}\label{weakker}
		Let $\lf$ be a~Leibniz algebra, and let $M$ be a~weak $\lf$-bimodule. Then $M\lf$ and $M^\lf$ are $\lf$-subbimodules of $M$.
	\end{lemma}
	From Lemma~\ref{weakker} we immediately obtain the following weak analogue of Theorem~3.14 in~\cite{F1}:
	\begin{proposition}\label{weakirred}
		Let $\lf$ be a~Leibniz algebra, and let $M$ be an irreducible weak
		$\lf$-bimodule. Then either $M$ is anti-symmetric or $M=M\lf$ and
		$M^\lf=0$.
	\end{proposition}
	Of course, Proposition~\ref{weakirred} is much weaker than its analogue for Leibniz bimodules. In the next example we obtain a~complete classification of the isomorphism classes of the irreducible weak Leibniz bimodules over an algebraically closed ground field. This is possible because the weak universal enveloping algebra of the one-dimensional Lie algebra is commutative, and thus there is a~geometric correspondence between irreducible modules and points in affine space.
	\begin{example}\label{weak1dim}
		{Let $\ef:=\fbb e$ be the one-dimensional Lie algebra over a~field
		$\fbb$. Observe that the weak universal enveloping algebra $\UL_{\weak}
		(\ef)$ is just a~polynomial algebra in two variables. Namely, if we set $x:=
		\ell_e$ and $y:=r_e$, then we have that}
		\[
			\UL_{\weak}(\ef)=\fbb\{x,y\}/\langle xy-yx\rangle=\fbb[x,y]\,,
		\]
		{where $\fbb\{x,y\}$ denotes the free associative $\fbb$-algebra in the variables $x$ and $y$, $\fbb[x,y]$ denotes the polynomial algebra in two commuting variables over $\fbb$, and $\langle xy-yx\rangle$ denotes the ideal generated by the element $xy-yx$.}

		{Now assume that $\fbb$ is algebraically closed. Then the weak $\ef$-bimodules correspond to (left) $\fbb[x,y]$-modules, and the isomorphism classes of the irreducible (left) $\fbb[x,y]$-modules are in bijection with the points of the affine space $\fbb^2$ via}
		\[
			(\alpha,\beta)\mapsto\fbb[x,y]/\langle x-\alpha,y-\beta\rangle\,.
		\]
		{Note that the latter are one-dimensional $\fbb[x,y]$-modules with underlying vector space $\fbb$, where $x$ acts via multiplication by $\alpha$ and $y$ acts via multiplication by $\beta$, i.e., $e$ acts on $\fbb$ from the left by multiplication with $\alpha$ and $e$ acts on $\fbb$ from the right by multiplication with $\beta$. We denote these weak $\ef$-bimodules by $_\alpha F_\beta$. In particular, we have shown that every irreducible weak $\ef$-bimodule over an algebraically closed ground field is one-dimensional.}

		{Note that $_\alpha F_\beta$ is symmetric exactly when $\alpha+\beta=0$, and it is anti-symmetric exactly when $\beta=0$. Consequently, the irreducible weak $\ef$-bimodule $_0 F_1$ is neither symmetric nor anti-symmetric. Moreover, for $M:=\,_0 F_1$ we have that $M_0 =M$ whose right action is clearly not trivial\footnote{The weak $\ef$-bimodule $M:=\,_0 F_1$ also shows that~\cite[Lemma~1.1]{FW1} does not hold for weak Leibniz bimodules.}}.

		{The following example shows that the tensor product of two ordinary
		Leibniz bimodules is not always an ordinary Leibniz bimodule. Namely, set
		$M:=\,_1 F_{-1}$ and $N:=\,_1 F_0$. Then $M\otimes N=\,_2 F_{-1}$ is irreducible, but neither symmetric nor anti-symmetric, and therefore $M\otimes N$ is only a weak $\ef$-bimodule (see Corollary~\ref{weakantisymmirred})}.
	\end{example}

	\begin{remark}\label{nonperfect}
		{Note that every non-perfect Leibniz algebra admits a~one-dimen\-sional
		Leibniz bimodule that is neither symmetric nor anti-symmetric. Namely, let $\lf$
		be a~non-perfect Leibniz algebra. Then $\lf/\lf\lf\ne 0$\footnote{Here $\lf\lf:=
		\langle xy\mid x,y\in\lf\rangle_\fbb$ is the \emph{derived subalgebra\/}
		of the Leibniz algebra $\lf$, and $\lf$ is called \emph{perfect\/} if $\lf=\lf\lf$.},
		and thus $(\lf/\lf\lf)^*\ne 0$. Now choose a~non-zero linear form on $\lf/\lf\lf$
		and lift it to a~non-zero linear form $\lambda$ on $\lf$. Then the one-dimensional weak $\lf$-bimodule $_0 F_\lambda$ on which $\lf$ acts trivially from the left and by $\xi\cdot x:=\lambda(x)\xi$ from the right is neither symmetric nor anti-symmetric.}

		{It should also be mentioned that it is not straightforward to generalize
		Example~\ref{weak1dim} to abelian Lie algebras of dimension greater than 1
		because already for a~two-dimensional abelian Lie algebra $\aaf:=\fbb e\oplus
		\fbb f$ the weak universal enveloping algebra is not commutative. Namely, we have that}
		\[
			\UL_{\weak}(\aaf)=\fbb\{x,y,u,v\}/(xy-yx,xu-ux,xv-vx,yu-uy,uv-vu)\,,
		\]
		{where $x:=\ell_e$, $y:=r_e$, $u:=\ell_f$,   $v:=r_f$, i.e., all variables commute except for $y$ and $v$.}
	\end{remark}
	It would be very interesting to classify all irreducible weak Leibniz bimodules up to isomorphism. We hope to come back to this problem at a~later time.

	Contrary to $\UL(\lf)$, which is only an augmented unital associative algebra, but similar to $\U(\lf_{\lie})$, $\UL_{\weak}(\lf)$ is a~cocommutive
	Hopf algebra:
	\begin{theorem}\label{hopf}
		Let $\lf$ be a~Leibniz algebra over a~field $\fbb$. Then $\UL_{\weak}(\lf)$
		is a~cocommutative Hopf algebra with comultiplication
		\[
			\Delta:\UL_{\weak}(\lf)\to\UL_{\weak}(\lf)\otimes\UL_{\weak}(\lf)
		\]
		defined by
		\[
			\Delta(\ell_x ):=\ell_x\otimes 1+1\otimes\ell_x
		\]
		and
		\[
			\Delta(r_x ):=r_x\otimes 1+1\otimes r_x\,,
		\]
		counit
		$\varepsilon:\UL_{\weak}(\lf)\to\fbb$ defined by
		\[
			\varepsilon(\ell_x ):=0=:\varepsilon(r_x )\,,
		\]
		and antipode $S:\UL_{\weak}(\lf)\to\UL_{\weak}(\lf)$
		defined by
		\[
			S(\ell_x ):=-\ell_x\quad\mbox{and}\quad S(r_x ):=-r_x\,.
		\]
	\end{theorem}

	\begin{proof}
		According to~\cite[Theorem III.2.4 and Example 3 on p.\ 56]{K}, the tensor algebra $T(\lf^\ell\oplus\lf^r )$ is a~cocommutative Hopf algebra with comultiplication, counit, and antipode defined as above. In order to enable us to lift these algebra homomorphisms to the factor algebra
		$\UL_{\weak}(\lf)$, we must show that $J(\lf)$ is a~coideal of $T(\lf^\ell
		\oplus\lf^r )$ (see~\cite[Definition~III.1.5]{K}) such that $S[J(\lf)]
		\subseteq J(\lf)$. For the counit by definition we have that $\varepsilon
		(\ell_x )=0=\varepsilon(r_x )$ for every $x\in\lf$. Moreover, we obtain that
		\begin{eqnarray*}
			\Delta(\ell_x\cdot r_y -r_y\cdot\ell_x -r_{xy}) & = & \Delta(\ell_x )\cdot
			\Delta(r_y )-\Delta(r_y )\cdot\Delta(\ell_x )-\Delta(r_{xy})\\
			& = & (\ell_x\otimes 1+1\otimes\ell_x )\cdot(r_y\otimes 1+1\otimes r_y )\\
			&& -(r_y\otimes 1+1\otimes r_y )\cdot(\ell_x\otimes 1+1\otimes\ell_x )\\
			&& -(r_{xy}\otimes 1+1\otimes r_{xy})\\
			& = & (\ell_x\cdot r_y )\otimes 1+\ell_x\otimes r_y +r_y\otimes\ell_x
			+1\otimes(\ell_x\cdot r_y )\\
			&& -(r_y\cdot\ell_x )\otimes 1-r_y\otimes\ell_x -\ell_x\otimes r_y -1
			\otimes(r_y\cdot\ell_x )\\
			&& -(r_{xy}\otimes 1+1\otimes r_{xy})\\
			& = & (\ell_x\cdot r_y -r_y\cdot\ell_x -r_{xy})\otimes 1+1\otimes
			(\ell_x\cdot r_y -r_y\cdot\ell_x - r_{xy})\\
			& \in & J(\lf)\otimes T(\lf^\ell\oplus\lf^r )+T(\lf^\ell\oplus\lf^r )\otimes J(\lf)
		\end{eqnarray*}
		for all elements $x,y\in\lf$. Similarly, we can show (as for universal enveloping algebras of Lie algebras) that
		\[
			\Delta(\ell_x\cdot\ell_y -\ell_y\cdot\ell_x -\ell_{xy})\in J(\lf)\otimes
			T(\lf^\ell\oplus\lf^r )+T(\lf^\ell\oplus\lf^r )\otimes J(\lf)
		\]
		for all elements $x,y\in\lf$.

		Finally, we have that
		\begin{eqnarray*}
			S(\ell_x\cdot r_y -r_y\cdot\ell_x -r_{xy}) & = & S(r_y )\cdot S(\ell_x )
			-S(\ell_x )\cdot S(r_y )-S(r_{xy})\\
			& = & (-r_y )\cdot(-\ell_x )-(-\ell_x )\cdot(-r_y )-(-r_{xy})\\
			& = & r_y\cdot\ell_x -\ell_x\cdot r_y +r_{xy}\in J(\lf)
		\end{eqnarray*}
		for all elements $x,y\in\lf$, and, similarly, we obtain that
		\[
			S(\ell_x\cdot\ell_y -\ell_y\cdot\ell_x -\ell_{xy})\in J(\lf)
		\]
		for all elements $x,y\in\lf$.
	\end{proof}

	\begin{remark}
		{As for the universal enveloping algebra of a~Lie algebra, the weak universal enveloping algebra of a~Leibniz algebra is generated by \emph{primitive elements\/}. Note that therefore it follows from~\cite[Proposition III.2.6 and (3.3)]{K} that the counit and antipode necessarily must be defined as in Theorem~\ref{hopf}.}

		{By virtue of the Cartier-Gabriel-Kostant theorem (see~\cite[Theorem 5.10.2]{EGNO}), the weak universal enveloping algebra $\UL_{\weak}(\lf)$ of a~Leibniz algebra $\lf$ over an algebraically closed field of characteristic zero is the universal enveloping algebra of the Lie algebra of its primitive elements, $\mathrm{Prim}
		(\UL_{\weak}(\lf))$. (Note that the unity is the only group-like element of $\UL_{\weak}(\lf)$
		because the latter is generated by primitive elements.)}
	\end{remark}
	It would be very useful to have an explicit description of
	$\mathrm{Prim}(\UL_{\weak}(\lf))$ in terms of a~Lie algebra naturally associated with $\lf$. Clearly, one such candidate would be the direct sum $\lf_{\lie}\oplus\lf_{\lie}$ of two copies of the canonical Lie algebra $\lf_{\lie}$ associated with $\lf$. But $\mathrm{Prim}(\UL_{\weak}(\lf))$ is not always isomorphic to $\lf_{\lie}\oplus\lf_{\lie}$ as the following example shows:
	\begin{example}
		{Let $\af=\fbb h\oplus\fbb e$ be the solvable left
		Leibniz algebra over a~field $\fbb$ with the multiplication
		$he=e$ and $hh=eh=ee=0$ (see~\cite[Example 2.3]{F1}). If we set $x:=\ell_h$, $y:=r_h$, and $z:=r_e$, then we have that}
		\[
			\UL_{\weak}(\af)=\fbb\{x,y,z\}/\langle xy-yx,[x,z]-z\rangle\,,
		\]
		{where $\fbb\{x,y,z\}$ denotes the free associative
		$\fbb$-algebra in the variables $x$, $y$, $z$, where $[x,z]:=xz-zx$
		is the commutator of $x$ and $z$, and where $\langle X\rangle$ denotes the two-sided ideal of $\fbb\{x,y,z\}$ that is generated by the set $X$.}

		{Then for degree reasons we have that $\{x+J(\af),y+J(\af), z+J(\af)\}$ is a~basis of the vector space $[\af^\ell\oplus\af^r
		+J(\af)]/J(\af)$. (Note that $he=e$, $eh=0$, and the relation
		\eqref{llm} imply $\ell_e =0$.) Consequently, we obtain that}
		\[
			\dim_\fbb(\af_{\lie}\oplus\af_{\lie})=2<\dim_\fbb[\af^\ell\oplus
			\af^r +J(\af)]/J(\af)\le\dim_\fbb\mathrm{Prim}(\UL_{\weak}(\af))\,,
		\]
		{which contradicts $\mathrm{Prim}(\UL_{\weak}(\af))\cong
		\af_{\lie}\oplus\af_{\lie}$.}
	\end{example}
	From the inclusion $J(\lf)\subseteq I(\lf)$ we obtain that there is an epimorphism
	\[
		\omega:\UL_{\weak}(\lf)\to\UL(\lf)
	\]
	of unital associative algebras. Let $\varepsilon_{\lie}:\U(\lf_{\lie})\to
	\fbb$ denote the counit of $\U(\lf_{\lie})$. Then we have the following composition of algebra homomorphisms:
	\[
		\UL_{\weak}(\lf) \stackrel{\omega}\to \UL(\lf) \stackrel{d_0}\to \U(\lf_{\lie}) \stackrel{\varepsilon_{\lie}}\to \fbb\,.
	\]
	In particular, $\varepsilon=\varepsilon_{\lie}\circ d_0\circ\omega$
	is the counit of $\UL_{\weak}(\lf)$. Similarly, we have
	\[
        \UL_{\weak}(\lf) \stackrel{\omega}\to \UL(\lf) \stackrel{d_1}\to \U(\lf_{\lie}) \stackrel{\varepsilon_{\lie}}\to \fbb
	\]
	and $\varepsilon=\varepsilon_{\lie}\circ d_1\circ\omega$ is the counit of $\UL(\lf)$ (see Proposition~\ref{augmented} for the definitions of the algebra homomorphisms $d_0$ and $d_1$).

	We obtain from Theorem~\ref{hopf} and~\cite[Proposition III.5.1]{K}
	that $\Mod_{\weak}^{\bi}(\lf)$ is a~symmetric monoidal category (see~\cite[Section~VII.1, pp.~162/163, Section~XI.1, pp.~252/253]{McL}, \cite[Definition XI.2.1]{K}, and~\cite[Definitions 8.1.1, 8.1.2, and
	8.1.12]{EGNO}):
	\begin{theorem}\label{monoidal}
		Let $\lf$ be a~Leibniz algebra. Then $\Mod_{\weak}^{\bi}(\lf)$ is a symmetric monoidal category.
	\end{theorem}
	In the following we will sketch a~more down-to-earth proof of
	Theorem~\ref{monoidal} that is neither using Theorem~\ref{hopf}
	nor~\cite[Proposition III.5.1]{K}, but instead relies on the well-known fact that the category $\vect$ of vector spaces is a symmetric monoidal category, i.e., we employ the associativity and commutativity constraints as well as the left and right unit constraints that make $\vect$ into a~symmetric monoidal category. In particular, then the pentagon, triangle, and hexagon axioms are clearly satisfied, and we only need to verify that the associativity and commutativity constraints as well as the left and right unit constraints of $\vect$ are homomorphisms of weak $\lf$-bimodules.

	Firstly, we show that for any three weak $\lf$-bimodules $L$, $M$, and $N$ the canonical isomorphism
	\[
		\alpha_{L,M,N}:(L\otimes M)\otimes N\to L\otimes(M\otimes N)\,,\,
		(l\otimes m)\otimes n\mapsto l\otimes(m\otimes n)
	\]
	is a~homomorphism of weak $\lf$-bimodules, i.e., $\alpha_{L,M,N}$
	is compatible with the left and right $\lf$-actions. Here we only verify the compatibility with the right $\lf$-action and leave the completely analogous proof for the compatibility with the left $\lf$-action to the interested reader. We have that
	\begin{align*}
		\alpha_{L,M,N}{ }&{ }([(l\otimes m)\otimes n]\cdot x) \\
        & = \alpha_{L,M,N}
		([(l\otimes m)\cdot x]\otimes n+(l\otimes m)\otimes(n\cdot x))\\
		& = \alpha_{L,M,N}([(l\cdot x)\otimes m]\otimes n+[l\otimes(m\cdot x)]
		\otimes n+(l\otimes m)\otimes(n\cdot x))\\
		& = (l\cdot x)\otimes(m\otimes n)+l\otimes[(m\cdot x)\otimes n]
		+l\otimes[m\otimes(n\cdot x)]\\
		& = (l\cdot x)\otimes(m\otimes n)+l\otimes[(m\otimes n)\cdot x]\\
		& = [l\otimes(m\otimes n)]\cdot x\\
		& = \alpha_{L,M,N}((l\otimes m)\otimes n)\cdot x
	\end{align*}
	for all elements $l\in L$, $m\in M$, $n\in N$, and $x\in\lf$.

	Secondly, we show that for any two weak $\lf$-bimodules $M$
	and $N$ the flip
	\[
		\gamma_{M,N}:M\otimes N\to N\otimes M\,,\, m\otimes n\mapsto n\otimes m
	\]
	is a~homomorphism of weak $\lf$-bimodules. Here we again only verify that $\gamma_{M,N}$ is compatible with the right $\lf$-action and leave the proof for the compatibility of $\gamma_{M,N}$ with the left $\lf$-action to the interested reader. We have that
	\begin{eqnarray*}
		\gamma_{M,N}((m\otimes n)\cdot x) & = & \gamma_{M,N}((m\cdot x)
		\otimes n+m\otimes(n\cdot x))\\
		& = & n\otimes(m\cdot x)+(n\cdot x)\otimes m\\
		& = & (n\cdot x)\otimes m+n\otimes(m\cdot x)\\
		& = & (n\otimes m)\cdot x\\
		& = & \gamma_{M,N}(m\otimes n)\cdot x
	\end{eqnarray*}
	for all elements $m\in M$, $n\in N$, and $x\in\lf$.

	Finally, we show that for every weak $\lf$-bimodule $M$
	the left unit
	\[
		\lambda_M:\fbb\otimes M\to M\,,\,1\otimes m\mapsto m
	\]
	and the right unit
	\[
		\rho_M:M\otimes\fbb\to M\,,\,m\otimes 1\mapsto m
	\]
	are homomorphisms of weak $\lf$-bimodules. Here we only show that $\lambda_M$ is compatible with the right $\lf$-action and leave the remaining proofs to the interested reader. Namely, let $x\in\lf$ and $m\in M$ be arbitrary. Then we have that
	\begin{eqnarray*}
		\lambda_M [(\beta\otimes m)\cdot x] & = & \lambda_M [\beta\otimes
		(m\cdot x)]=\beta(m\cdot x)\\
		& = & (\beta m)\cdot x=\lambda_M (\beta\otimes m)\cdot x
	\end{eqnarray*}
	for all elements $\beta\in\fbb$, $m\in M$, and $x\in\lf$.

	As to be expected, the tensor product of weak Leibniz bimodules is compatible with direct sums:
	\begin{proposition}\label{distrib}
		Let $\lf$ be a~Leibniz algebra, and let $L$, $M$, $N$ be weak
		$\lf$-bimodules. Then there are natural isomorphisms
		\[
			L\otimes(M\oplus N)\cong(L\otimes M)\oplus(L\otimes N)
		\]
		and
		\[
			(L\oplus M)\otimes N\cong(L\otimes N)\oplus(M\otimes N)
		\]
		of weak $\lf$-bimodules.
	\end{proposition}
	As for Theorem~\ref{monoidal}, the natural isomorphisms are again the ones in the symmetric monoidal category of $\vect$, and similar to the previous paragraph, one can verify that these are morphisms in $\Mod_{\weak}^{\bi}(\lf)$.

	Clearly, the category $\Mod^{\bi}(\lf)$ of $\lf$-bimodules is a full subcategory of the category $\Mod_{\weak}^{\bi}(\lf)$ of weak
	$\lf$-bimodules. In fact, more is true: The forgetful (or underlying or inclusion or restriction) functor
	\[
		U:\Mod^{\bi}(\lf)\to\Mod_{\weak}^{\bi}(\lf)
	\]
	admits the left adjoint induction functor
	\[
		I:\Mod_{\weak}^{\bi}(\lf)\to\Mod^{\bi}(\lf)
	\]
	defined by $M\mapsto \UL(\lf)\otimes_{\UL_{\weak}(\lf)}M$, where
	$\UL(\lf)$ is a~right $\UL_{\weak}(\lf)$-module via the epimorphism
	$\omega:\UL_{\weak}(\lf)\to\UL(\lf)$ of unital associative algebras. This is indeed an adjunction:
	\[
		\Hom_{\Mod^{\bi}(\lf)}(I(M), N)\cong\Hom_{\Mod_{\weak}^{\bi}(\lf)}
		(M,U(N))\,.
	\]
	In fact, this natural isomorphism is the usual change-of-rings adjunction
	(Frobenius reciprocity or Shapiro's lemma).

	Recall that a~subcategory is called \textbf{reflective} if the inclusion functor admits a~left adjoint (see~\cite[p.\ 91]{McL}). Thus, we have
	\begin{proposition}\label{reflectiv}
		Let $\lf$ be a~Leibniz algebra. Then $\Mod^{\bi}(\lf)$ is a~reflective subcategory of $\Mod_{\weak}^{\bi}(\lf)$.
	\end{proposition}
	Clearly, the inclusion functor $U:\Mod^{\bi}(\lf)\to\Mod_{\weak}^{\bi}(\lf)$
	is fully faithful. Hence, we obtain from~\cite[Theorem 1, p.\ 90]{McL}
	the following result:
	\begin{proposition}
		If $\lf$ is a~Leibniz algebra, then for every $\lf$-bimodule $M$ we have the natural isomorphism
		\[
			\UL(\lf)\otimes_{\UL_{\weak}(\lf)}U(M)\cong M
		\]
		of $\lf$-bimodules.
	\end{proposition}
	It would be interesting to know whether the category of weak Leibniz bimodules is the \textit{smallest} symmetric monoidal category containing the category of Leibniz bimodules as a~reflective subcategory.

	Let $\fdmod_{\weak}^{\bi}(\lf)$ denote the category of finite-dimensional weak $\lf$-bimodules for a~Leibniz algebra $\lf$ which is a~full subcategory of the category $\Mod_{\weak}^{\bi}(\lf)$ of all weak
	$\lf$-bimodules. If $\UL_{\weak}(\lf)-\fdmod$ denotes the category of finite-dimensional unital left $\UL_{\weak}(\lf)$-modules, then we have the following equivalence of categories:
	\[
		\fdmod_{\weak}^{\bi}(\lf)\cong\UL_{\weak}(\lf)-\fdmod\,.
	\]
	As a~consequence, we obtain that the category of finite-dimensional weak Leibniz bimodules is a~locally finite $\fbb$-linear abelian category (see~\cite[Chapter VIII]{McL} or~\cite[Definitions~1.8.1, 1.2.2, and 1.3.1, respectively]{EGNO}):

	\begin{proposition}\label{locfinab}
		Let $\lf$ be a~Leibniz algebra over a~field $\fbb$. Then
		$\fdmod_{\weak}^{\bi}(\lf)$ is a~locally finite $\fbb$-linear abelian category.
	\end{proposition}

	\begin{remark}
		It follows from Example~\ref{weak1dim} that $\fdmod_{\weak}^{\bi} (\lf)$ is not always finite (see condition (iv) in~\cite[Definition 1.8.6]{EGNO}).
	\end{remark}

	Moreover, it is clear that $\fdmod_{\weak}^{\bi}(\lf)$ is a~full monoidal subcategory of the monoidal category $\Mod_{\weak}^{\bi}(\lf)$ (see~\cite[Definition 2.1.4]{EGNO}\footnote{Note that here we use~\cite[Definition 2.2.8]{EGNO} (or~\cite[pp.~162/163]{McL} resp. \cite[Definition XI.2.1]{K}) instead of~\cite[Definition 2.1.1]{EGNO} for the definition of a~monoidal category. Then a~\emph{full monoidal subcategory} of a~monoidal category $(\ccal, \otimes, 1,\alpha, \lambda, \rho)$ is a~sextuple $(\dcal, \otimes, 1, \alpha, \lambda, \rho)$ such that $\dcal$ is a~full subcategory of $\ccal$ that contains the unit $1$ and is closed under tensor products of objects and morphisms.}). We conclude this section by showing that $\fdmod_{\weak}^{\bi}(\lf)$ is rigid, i.e., every finite-dimensional weak $\lf$-bimodule has a~left and a right dual (see~\cite[Definitions 2.10.1, 2.10.2, and 2.10.11]{EGNO}).

	First, we need to show that the linear dual of a~weak Leibniz bimodule is again a~weak Leibniz bimodule. More generally, we prove that the vector space of linear transformations between weak Leibniz bimodules is a~weak Leibniz bimodule.

	Let $\lf$ be a~left Leibniz algebra, and let $M$ and $N$ be $\lf$-bimodules. We define left and right actions of $\lf$ on the vector space $\Hom_\fbb
	(M,N)$ of linear transformations from $M$ to $N$ as follows:
	\[
		(x\cdot f)(m):=x\cdot f(m)-f(x\cdot m)
	\]
	and
	\[
		(f\cdot x)(m):=f(m)\cdot x-f(m\cdot x)
	\]
	for all elements $x\in\lf$, $f\in\Hom_\fbb(M,N)$, and $m\in M$. Then we have the following result (cf.\ Proposition~\ref{tensprod}):
	\begin{proposition}\label{hom}
		Let $\lf$ be a~left Leibniz algebra over a~field $\fbb$, and let $M$ and
		$N$ be $\lf$-bimodules. Then $\Hom_\fbb(M,N)$ is a~weak $\lf$-bimodule. In particular, if $M$ and $N$ are both symmetric or both anti-symmetric, then $\Hom_\fbb(M,N)$ is a~symmetric or anti-symmetric
		$\lf$-bimodule, respectively.
	\end{proposition}

	\begin{proof}
		As is true for Lie modules, $\Hom_\fbb(M,N)$ satisfies \textrm{(LLM)} (see also the first part of the proof of~\cite[Lemma 1.4\,(b)]{FW1} for the special case $M=\lf_{\ad}$).

		Next, let us verify that (LML) holds. Namely, as (LML) is satisfied for $M$
		and $N$, we have that
		\begin{eqnarray*}
			[f\cdot(xy)]
			(m) & = & f(m)\cdot(xy)-f[m\cdot(xy)]\\
			& = & x\cdot[f(m)\cdot y]-[x\cdot f(m)]\cdot y-f[x\cdot(m\cdot y)]+f[(x\cdot m)\cdot y]\\
			& = & x\cdot[f(m)\cdot y]-x\cdot f(m\cdot y)-f(x\cdot m)\cdot y+f[(x\cdot m)\cdot y]\\
			&& -[x\cdot f(m)]\cdot y+f(x\cdot m)\cdot y+x\cdot f(m\cdot y)-f[x\cdot(m\cdot y)]\\
			& = & x\cdot(f\cdot y)(m)-(f\cdot y)(x\cdot m)-(x\cdot f)(m)\cdot y+(x\cdot f)(m\cdot y)\\
			& = & [x\cdot(f\cdot y)](m)-[(x\cdot f)\cdot y](m)
		\end{eqnarray*}
		for all elements $f\in\Hom_\fbb(M,N)$, $x,y\in\lf$, and $m\in M$.

		Finally, let us consider (ZD). Because (ZD) holds for $N$, we obtain that
		\begin{eqnarray*}
			[(x\cdot f+f\cdot x)\cdot y]
			(m) & = & (x\cdot f+f\cdot x)(m)\cdot y
			-(x\cdot f+f\cdot x)(m\cdot y)\\
			& = & (x\cdot f)(m)\cdot y+(f\cdot x)(m)\cdot y-(x\cdot f)(m\cdot y)-(f\cdot x)(m\cdot y)\\
			& = & [x\cdot f(m)]\cdot y-f(x\cdot m)\cdot y+[f(m)\cdot x]\cdot y-f(m\cdot x)\cdot y\\
			&& -x\cdot f(m\cdot y)+f[x\cdot(m\cdot y)]-f(m\cdot y)\cdot x+f[(m\cdot y)\cdot x]\\
			& = & \underbrace{[x\cdot f(m)+f(m)\cdot x]\cdot y}_{=0}-f(x\cdot m+m\cdot x)\cdot y\\
			&& -x\cdot f(m\cdot y)-f(m\cdot y)\cdot x+f[x\cdot(m\cdot y)+(m\cdot y)\cdot x]\\
			& = & f[x\cdot(m\cdot y)+(m\cdot y)\cdot x]-[x\cdot f(m\cdot y)+f(m\cdot y)\cdot x]\\
			&& -f(x\cdot m+m\cdot x)\cdot y
		\end{eqnarray*}
		for all elements $f\in\Hom_\fbb(M,N)$, $x,y\in\lf$, and $m\in M$. From the last expression we conclude that $\Hom_\fbb(M,N)$ is an (anti-)symmetric
		$\lf$-bimodule if both $M$ and $N$ are (anti-)symmetric.
	\end{proof}

	\begin{remark}
		It should be noted that in the above proof of (LML), we only used that
		(LML) is satisfied for each of the bimodules. Similarly, the proof of (LLM) only uses that (LLM) is satisfied for each of the bimodules. Moreover, it is noteworthy that in the last part of the proof of Proposition~\ref{hom} only the second bimodule needs to satisfy (ZD).
	\end{remark}

	The proof of Proposition~\ref{hom} shows that the vector space of linear transformations between weak Leibniz bimodules is again a~weak Leibniz bimodule:

    \begin{proposition}
		Let $\lf$ be a~Leibniz algebra over a~field $\fbb$. If $M$ and $N$ are weak
		$\lf$-bimodules, then $\Hom_\fbb(M,N)$ is a~weak $\lf$-bimodule. In particular, the linear dual $M^*$ of a~weak $\lf$-bimodule is a~weak $\lf$-bimodule.
	\end{proposition}

	Let $\lf$ be a~Leibniz algebra over a~field $\fbb$, and let $M$ be a~weak
	$\lf$-bimodule. Define the contractions
	\[
		\ev_M:M^*\otimes M\to\fbb\,,\,\mu\otimes m\mapsto\mu(m)
	\]
	and
	\[
		\ev_M ':M\otimes M^*\to\fbb\,,\,m\otimes\mu\mapsto\mu(m)\,.
	\]
	If $M$ is finite dimensional, choose a~basis $\{m_1,\dots,m_d\}$ of
	$M$ and the canonical dual basis $\{m_1^*,\dots,m_d^*\}$ such that $m_i^* (m_j )=\delta_{ij}$ for all $i,j\in\{1,\dots,d\}$, where
	$\delta_{ij}$ is the usual Kronecker delta. Then we can define the embeddings
	\[
		\coev_M:\fbb\to M\otimes M^*\,,\,1\mapsto\sum_{i=1}^d m_i
		\otimes m_i^*
	\]
	and
	\[
		\coev_M ':\fbb\to M^*\otimes M\,,\,1\mapsto\sum_{i=1}^d m_i^*
		\otimes m_i
	\]
	of $\fbb$-vector spaces. Note that
    \[
        \ev_M\circ\coev_M '
        = (\dim_\fbb M)\cdot\id_\fbb
        \quad \text{and} \quad
        \ev_M '\circ\coev_M
        = (\dim_\fbb M)\cdot\id_\fbb.
    \]
    In particular, $\ev_M$, $\ev_M '$ are surjective and $\coev_M$, $\coev_M '$ are injective.
	\begin{theorem}\label{rigid}
		Let $\lf$ be a~Leibniz algebra. Then $\fdmod_{\weak}^{\bi}(\lf)$ is a~rigid symmetric monoidal category.
	\end{theorem}

	\begin{proof}
		It is well known that the category of finite-dimensional vector spaces is a rigid symmetric monoidal category with the contractions and embeddings defined as above (see~\cite[Example 2.10.12]{EGNO}), where the linear dual $V^*$ of a~finite-dimensional vector space $V$ is a~left and right dual of $V$. Consequently, it is enough to prove that $\ev_M$, $\ev_M '$,
		$\coev_M$, and $\coev_M '$ are homomorphisms of weak $\lf$-bimodules for every finite-dimensional weak $\lf$-bimodule $M$.

		Let $M$ be a~finite-dimensional weak $\lf$-bimodule. We will show that
		$\ev_M$ and $\coev_M '$ are compatible with the right action of $\lf$ and leave the remaining proofs, which are very similar, to the interested reader.

		Firstly, we prove that $\ev_M$ is compatible with the right action of $\lf$. Since $\lf$ acts trivially on $\fbb$, we need to show that $\ev_M [(\mu
		\otimes m)\cdot x]=0$ for all elements $x\in\lf$, $\mu\in M^*$, and
		$m\in M$:
		\begin{eqnarray*}
			\ev_M [(\mu\otimes m)\cdot x] & = & \ev_M [(\mu\cdot x)\otimes m+
			\mu\otimes(m\cdot x)]\\
			& = & (\mu\cdot x)(m)+\mu(m\cdot x)\\
			& = & -\mu(m\cdot x)+\mu(m\cdot x)=0\,.
		\end{eqnarray*}
		Next, we prove that $\coev_M'$ is compatible with the right action of $\lf$. For this we first need to show that if $m_i\cdot x=\sum_{j=1}^d\xi_{ij}
		m_j$, then $m_i^*\cdot x=-\sum_{j=1}^d\xi_{ji}m_j^*$ for every integer $i\in\{1,\dots,d\}$. Namely, let $m=\sum_{k=1}^d\beta_k m_k$ be arbitrary. Then we have that
		\begin{eqnarray*}
			(m_i^*\cdot x)(m) & = & -m_i^* (m\cdot x)\\
			& = & -\sum_{k=1}^d\beta_k m_i^* (m_k\cdot x)\\
			& = & -\sum_{k,l=1}^d\beta_k\xi_{kl}m_i^* (m_l )\\
			& = & -\sum_{k=1}^d\beta_k\xi_{ki}\\
			& = & -\sum_{j,k=1}^d\beta_k\xi_{ji}m_j^* (m_k )\\
			& = & -\sum_{j=1}^d\xi_{ji}m_j^* (m)\\
			& = & -\left(\sum_{j=1}^d\xi_{ji}m_j^*\right)(m)\,,
		\end{eqnarray*}
		and therefore we obtain that
		\begin{eqnarray*}
			\sum_{i=1}^d m_i^*\otimes(m_i\cdot x) & = & \sum_{i=1}^d
			\left[m_i^*\otimes\left(\sum_{j=1}^d\xi_{ij}m_j\right)\right]\\
			& = & \sum_{j=1}^d\left[\left(\sum_{i=1}^d\xi_{ij}m_i^*\right)
			\otimes m_j\right]\\
			& = & -\sum_{j=1}^d (m_j^*\cdot x)\otimes m_j\,.
		\end{eqnarray*}
		Now we can prove that $\coev_M '$ is compatible with the right action of $\lf$. Since $\lf$ acts trivially on $\fbb$, we need to show that $\coev_M '(\sigma)
		\cdot x=0$ for all elements $x\in\lf$ and $\sigma\in\fbb$:
		\begin{eqnarray*}
			\coev_M '(\sigma)\cdot x & = & \sigma\left(\sum_{i=1}^d m_i^*\otimes m_i\right)\cdot x\\
			& = & \sigma\sum_{i=1}^d [(m_i^*\cdot x)\otimes m_i +m_i^*\otimes(m_i\cdot x)]=0\,,
		\end{eqnarray*}
		as desired.
	\end{proof}

	\begin{remark}\label{pivotal}
		{One can also show for every finite-dimensional weak $\lf$-bimodule
		$M$ that the natural isomorphism $M\cong M^{**}$ of vector spaces is compatible with the left and right $\lf$-actions, i.e., $\fdmod_{\weak}^{\bi}
		(\lf)$ is \emph{pivotal\/} (see~\cite[Definition 4.7.8]{EGNO}).}
	\end{remark}
	According to Proposition~\ref{locfinab} and Theorem~\ref{rigid}, for every
	Leibniz algebra $\lf$ over a~field $\fbb$ we have that $\fdmod_{\weak}^{\bi}
	(\lf)$ is a~locally finite $\fbb$-linear abelian rigid symmetric monoidal category. Moreover, clearly, the bifunctor
	\[
		\otimes:\fdmod_{\weak}^{\bi}(\lf)\times\fdmod_{\weak}^{\bi}(\lf)
		\to\fdmod_{\weak}^{\bi}(\lf)
	\]
	is $\fbb$-bilinear on morphisms and $\End_\lf(\fbb)\cong\fbb$. Hence, we obtain that $\fdmod_{\weak}^{\bi}(\lf)$ is a~\emph{tensor category} over $\fbb$ in the sense of Etingof et al.\footnote{Note that in Kassel's book a~tensor category is the same as a~monoidal category (see~\cite[Definition~XI.2.1]{K}).}
	(see~\cite[Definition 4.1.1]{EGNO}), and therefore it follows from~\cite[Proposition 4.2.1]{EGNO} that the bifunctor $\otimes$
	is biexact\footnote{A bifunctor $\ccal\times\ccal\to\ccal$ is
	\emph{biexact\/} if it is exact in both factors.}. Consequently,
	$\fdmod_{\weak}^{\bi}(\lf)$ is a~\emph{ring category\/} over
	$\fbb$ (see~\cite[Definition 4.2.3]{EGNO}):
	\begin{theorem}\label{ring}
		Let $\lf$ be a~Leibniz algebra over a~field $\fbb$. Then
		$\fdmod_{\weak}^{\bi}(\lf)$ is a~ring category over $\fbb$.
	\end{theorem}
	In Section~\ref{grothendieck} we will study the Grothendieck ring of
	$\fdmod_{\weak}^{\bi}(\lf)$ and a~certain non-associative Grothendieck ring defined for the category $\fdmod^{\bi}(\lf)$ of all finite-dimensional
	$\lf$-bimodules and compare them with the Grothendieck ring of the category $\fdmod(\lf_{\lie})$ of finite-dimensional $\lf_{\lie}$-modules. Before we can do this, we need to introduce a~tensor product on the category of (finite-dimensional) Leibniz bimodules.

	\section{Truncated tensor products for Leibniz bimodules}\label{truncated}

	In this section we define for $\lf$-bimodules $M$ and $N$
	over a~left Leibniz algebra $\lf$ two truncated tensor products
	$M\ovo N$ and $M\und N$ which both are again $\lf$-bimodules. According to \eqref{*} in Proposition~\ref{tensprod}, we must ensure that
	\[
		(x\cdot m+m\cdot x)\otimes (n\cdot y)+(m\cdot y)\otimes
		(x\cdot n+n\cdot x)=0
	\]
	holds for all elements $x,y\in\lf$, $m\in M$, and $n\in N$.

    Let $T(M,N)$ be the smallest subspace of the tensor product 	$M\otimes N$ that contains
	\[
		S(M,N):=\{(x\cdot m+m\cdot x)\otimes (n\cdot y)+(m\cdot y)
		\otimes(x\cdot n+n\cdot x)\mid x,y\in\lf;m\in M,n\in N\}
	\]
	and that is closed under the left and right $\lf$-actions as defined in Section~\ref{weak}, i.e., $T(M,N)$ is the weak $\lf$-bimodule of $M\otimes N$ generated by $S(M,N)$. Now define
	\[
		M\ovo N:=(M\otimes N)/T(M,N)\,.
	\]

	On the other hand, set
	\[
		T_0 (M,N):=M_0\otimes N\lf+M\lf\otimes N_0\,,
	\]
	where $M_0$ and $N_0$ are the \emph{anti-symmetric kernels\/}
	of $M$ and $N$, respectively, which both are anti-symmetric
	$\lf$-subbimodules (see~\cite[Proposition 3.12]{F1}), and $M\lf$
	is the subspace of $M$ that is spanned by $\{m\cdot x\mid m\in
	M,x\in\lf\}$, and similarly, for $N\lf$. Since every factor in the tensor products of $T_0 (M,N)$ is an $\lf$-subbimodule, we obtain that $T_0 (M,N)$ is a~weak
	$\lf$-subbimodule of $M\otimes N$, and we define
	\[
		M\und N:=(M\otimes N)/T_0 (M,N)\,.
	\]

	\begin{remark}
		{Note that the truncated tensor product $M\ovo N$ even makes sense for weak $\lf$-bimodules $M$ and $N$. On the other hand, this is not clear for the truncated tensor product
		$M\und N$, although $M\lf$ and $N\lf$ are weak $\lf$-subbimodules of weak $\lf$-bimodules $M$ and $N$, respectively (see Lemma~\ref{weakker}), but this seems not to be the case for the anti-symmetric kernels $M_0$ and $N_0$, respectively (see also our discussion after Proposition~\ref{ab}).}
	\end{remark}
	It follows from Proposition~\ref{tensprod} that the truncated tensor products of Leibniz bimodules are again Leibniz bimodules:
	\begin{proposition}
		Let $\lf$ be a~Leibniz algebra, and let $M$ and $N$ be
		$\lf$-bimodules. Then the truncated tensor products $M\ovo N$
		and $M\und N$ are $\lf$-bimodules.
	\end{proposition}

	\begin{remark}
		{Observe that factoring out $T_0 (M,N)=M_0\otimes
		N\lf+M\lf\otimes N_0$ is factoring out more than $T(M,N)$, because in the definition of $S(M,N)$ the two summands are linked, while in $T_0 (M,N)$ these terms are not linked. Consequently, we have the inclusion $T(M,N)\subseteq
		T_0 (M,N)$, but since we could not find any example where
		$T(M,N)$ is properly contained in $T_0 (M,N)$, at the moment it is not clear to us whether both truncated tensor products always coincide.}
	\end{remark}
	It is not surprising that both truncated tensor products are commutative:
	\begin{proposition}\label{comm}
		Let $\lf$ be a~Leibniz algebra, and let $M$ and $N$ be
		$\lf$-bimodules. Then there are natural isomorphisms $M
		\ovo N\cong N\ovo M$ and $M\und N\cong N\und M$ of
		$\lf$-bimodules.
	\end{proposition}

	\begin{proof}
		Let $\gamma:M\otimes N\to N\otimes M$ be the natural isomorphism of weak $\lf$-bimodules defined by $\gamma(m\otimes n):=n\otimes m$ (cf.\ the down-to-earth proof of Theorem~\ref{monoidal}). It is easy to see that $\gamma[T(M,N)]=T(N,M)$. Now consider the homomorphism $\overline{\gamma}:=\overline{\eta}\circ
		\gamma:M\otimes N\to N\ovo M$ of weak $\lf$-bimodules, where
		$\overline{\eta}:N\otimes M\to N\ovo M$ is the natural epimorphism. As a~composition of surjective functions, $\overline{\gamma}$
		is also surjective. So it remains to prove that $\Ker(\overline{\gamma})
		=T(M,N)$ in order to establish the first assertion. But the latter follows from $\gamma[T(M,N)] =T(N,M)$.

		Similar to the above proof, we observe that $\gamma[T_0 (M,N)]
		=T_0 (N,M)$, and then we consider the homomorphism
		$\underline{\gamma}:=\underline{\eta}\circ\gamma:M\otimes N
		\to N\und M$ of weak $\lf$-bimodules, where $\underline{\eta}:
		N\otimes M\to N\und M$ is the natural epimorphism. Of course,
		$\underline{\gamma}$ is surjective. So it remains to prove that
		$\Ker(\underline{\gamma})=T_0 (M,N)$ which immediately follows from the fact that $\gamma[T_0 (M,N)]=T_0 (N,M)$.
	\end{proof}

	It should be noted that neither of the truncated tensor products is associative. We will come back to this later in this section
	(see Proposition~\ref{nonasso}) and in Section~\ref{grothendieck}
	(see Corollary~\ref{nonasso2}). But both truncated tensor products are compatible with direct sums:
	\begin{proposition}
		Let $\lf$ be a~Leibniz algebra, and let $L$, $M$, and $N$ be
		$\lf$-bimodules. Then there are natural isomorphisms
		\begin{gather*}
			L\ovo(M\oplus N)\cong(L\ovo M)\oplus(L\ovo N)\,, \\
			L\und(M\oplus N)\cong(L\und M)\oplus(L\und N)\,, \\
			(L\oplus M)\ovo N\cong(L\ovo N)\oplus(M\ovo N)\,, \\
			(L\oplus M)\und N\cong(L\und N)\oplus(M\und N)
		\end{gather*}
		of $\lf$-bimodules.
	\end{proposition}

	\begin{proof}
		According to Proposition~\ref{comm}, it is enough to prove the first two statements. By virtue of Proposition~\ref{distrib}, we have the natural isomorphism
		\[
			\delta:L\otimes(M\oplus N)\to(L\otimes M)\oplus(L\otimes N)\,,\, l\otimes(m,n)\mapsto(l\otimes m,l\otimes n)
		\]
		of weak $\lf$-bimodules.

		In order to establish the first isomorphism, it remains to show that
		$\delta$ lifts to the corresponding truncated tensor products. For this, we will need that
		\[
			\delta[S(L,M\oplus N)]\subseteq S(L,M)\times S(L,N)\,,
		\]
		which in turn implies that
		\[
			\delta[T(L,M\oplus N)]\subseteq T(L,M)\times T(L,N)\,.
		\]
		Every element $s\in S(L,M\oplus N)$ can be written as
		\[
			s=(x\cdot l+l\cdot x)\otimes(m,n)\cdot y+l\cdot y\otimes
			[x\cdot(m,n)+(m,n)\cdot x]
		\]
		for some elements $l\in L$, $m\in M$, $n\in N$, and $x,y\in\lf$. Then we have that
		\begin{eqnarray*}
			\delta[(x\cdot l+l\cdot x)\otimes(m,n)\cdot y] & = & \delta
			(x\cdot l\otimes(m\cdot y,n\cdot y)+l\cdot x\otimes(m\cdot y,n\cdot y))\\
			& = & (x\cdot l\otimes m\cdot y,x\cdot l\otimes n\cdot y)+(l\cdot x\otimes m
			\cdot y,l\cdot x\otimes n\cdot y)\\
			& = & ((x\cdot l+l\cdot x)\otimes m\cdot y,(x\cdot l+l\cdot x)\otimes n\cdot y)
		\end{eqnarray*}
		and
		\begin{align*}
			\delta { }&{ } (l\cdot y\otimes[x\cdot(m,n)+(m,n)\cdot x]) \\
            & =\delta
			(l\cdot y\otimes(x\cdot m,x\cdot n)+l\cdot y\otimes(m\cdot x,n\cdot x))\\
			& = (l\cdot y\otimes m\cdot x,l\cdot y\otimes n\cdot x)+(l\cdot y\otimes x\cdot m,l\cdot y\otimes x\cdot n)\\
			& = ((l\cdot y\otimes(x\cdot m+m\cdot x),l\cdot y\otimes(x\cdot n+n\cdot x))\,,
		\end{align*}
		which implies that
		\begin{eqnarray*}
			\delta(s) & = &
			((x\cdot l+l\cdot x)\otimes m\cdot y+l\cdot y\otimes(x\cdot m+m\cdot x),\\
			&&(x\cdot l+l\cdot x)\otimes n\cdot y+l\cdot y\otimes(x\cdot n+n\cdot x))\\
			& \in & S(L,M)\times S(L,N)\,.
		\end{eqnarray*}
		It follows directly from the additivity of the functors $M\mapsto M_0$
		and $M\mapsto M\lf$ that $\delta$ lifts to the second isomorphism.
	\end{proof}

	In the remainder of this section we will discuss when the two truncated tensor products coincide or are non-zero. In particular, the former happens when one of the bimodules is symmetric or anti-symmetric. Recall that an $\lf$-bimodule $M$ is symmetric exactly when $M_0 =0$, and $M$ is anti-symmetric exactly when
	$M\lf=0$.
	\begin{theorem}\label{main}
		Let $\lf$ be a~Leibniz algebra, and let $M$ and $N$ be $\lf$-bimodules. Then the following statements hold:
		\begin{enumerate}
			\item[\textrm{(a)}] If $M$ is symmetric, then
			\[
				M\ovo N=M\und N=(M\otimes N)/(\lf M\otimes N_0 )\,.
			\]
			\item[\textrm{(b)}] If $M$ is anti-symmetric, then
			\[
				M\ovo N=M\und N=(M\otimes N)/(\lf M\otimes N\lf)\,.
			\]
			\item[\textrm{(c)}] If $N$ is symmetric, then
			\[
				M\ovo N=M\und N=(M\otimes N)/(M_0\otimes\lf N)\,.
			\]
			\item[\textrm{(d)}] If $N$ is anti-symmetric, then
			\[
				M\ovo N=M\und N=(M\otimes N)/(M\lf\otimes\lf N)\,.
			\]
		\end{enumerate}
	\end{theorem}

	\begin{proof}
		Note that $X\lf=\lf X$ for a~symmetric $\lf$-bimodule $X$ and $Y_0
		=\lf Y$ for an anti-symmetric $\lf$-bimodule $Y$. This in conjunction with the definition of $T_0 (M,N)$ yields the statements (a) -- (d) for
		$M\und N$. By virtue of Proposition~\ref{comm}, it remains to prove that $T(M,N)=T_0 (M,N)$ if $M$ is either symmetric or anti-symmetric.

		Suppose that $M$ is symmetric. Then $M_0 =0$, and thus $T_0
		(M,N)=M\lf\otimes N_0$ and
		\[
			S(M,N)=\{(m\cdot y)\otimes(x\cdot n+n\cdot x)\mid x,y\in\lf;
			m\in M,n\in N\}\,.
		\]
		By using (LML) and the argument in the proof of Proposition 3.12
		in~\cite{F1}, we see then that $\langle S(M,N)\rangle_\fbb$ is a~weak
		$\lf$-subbimodule of $M\otimes N$, and therefore we obtain that
		\[
            T(M,N)
            =\langle S(M,N)\rangle_\fbb
            =M\lf\otimes N_0
            =T_0 (M,N).
        \]
        The proof in the other case works very similar and is left to the interested reader.
	\end{proof}

	In particular, we have for two (anti-)symmetric Leibniz bimodules that both truncated tensor products coincide with the ``natural'' tensor product defined in Section~\ref{weak}:
	\begin{corollary}\label{sym}
		Let $\lf$ be a~Leibniz algebra, and let $M$ and $N$ be $\lf$-bimodules. Then the following statements hold:
		\begin{enumerate}
			\item[\textrm{(a)}] If $M$ and $N$ are symmetric, then $M\ovo N
			=M\und N=M\otimes N$.
			\item[\textrm{(b)}] If $M$ and $N$ are anti-symmetric, then $M
			\ovo N=M\und N=M\otimes N$.
			\item[\textrm{(c)}] If $M$ is symmetric and $N$ is anti-symmetric, then
			\[
				M\ovo N=M\und N=(M\otimes N)/(\lf M
				\otimes\lf N)\,.
			\]
			\item[\textrm{(d)}] If $M$ is anti-symmetric and $N$ is symmetric, then
			\[
				M\ovo N=M\und N=(M\otimes N)/(\lf M
				\otimes\lf N)\,.
			\]
		\end{enumerate}
	\end{corollary}

	\begin{proof}
		(a) is an immediate consequence of Theorem~\ref{main}\,(a) or
		(c), and similarly, (b) follows from Theorem~\ref{main}\,(b) or
		(d). Moreover, (c) is a~special case of Theorem~\ref{main}\,(a)
		or (d), and (d) is a~special case of Theorem~\ref{main}\,(b) or
		(c).
	\end{proof}

	Recall that a~\emph{trivial Leibniz bimodule\/} is a~Leibniz bimodule with trivial left and right actions. In particular, the ground field
	$\fbb$ of a~Leibniz algebra $\lf$ with trivial left and right $\lf$-actions is called the \emph{one-dimensional trivial $\lf$-bimodule\/} and will be denoted by $F_0^{s/a}$. Note that trivial Leibniz bimodules are symmetric and anti-symmetric. As an immediate consequence of
	Theorem~\ref{main} we obtain that the truncated tensor product of a~trivial bimodule with an arbitrary bimodule coincides with the
	``natural'' tensor product defined in Section~\ref{weak}:
	\begin{corollary}\label{triv}
		Let $\lf$ be a~Leibniz algebra over a~field $\fbb$, and let $M$ and
		$N$ be $\lf$-bimodules. If one of the $\lf$-bimodules is trivial, then $M\ovo N=M\und N=M\otimes N$. In particular,
        \[
            M\ovo
    		F_0^{s/a}=M\und F_0^{s/a}=M\otimes F_0^{s/a}\cong M\cong
    		F_0^{s/a}\otimes M=F_0^{s/a}\ovo M=F_0^{s/a}\und M.
        \]
	\end{corollary}
	It is not necessarily true that the truncated tensor products are non-zero:
	\begin{corollary}\label{irred}
		Let $\lf$ be a~Leibniz algebra, and let $M$ and $N$ be non-trivial irreducible $\lf$-bimodules. Then the following statements hold:
		\begin{enumerate}
			\item[\textrm{(a)}] If $M$ is symmetric and $N$ is anti-symmetric, then $M\ovo N=M\und N=0$.
			\item[\textrm{(b)}] If $M$ is anti-symmetric and $N$ is symmetric, then $M\ovo N=M\und N=0$.
		\end{enumerate}
	\end{corollary}

	\begin{proof}
		Observe that under either of the two hypotheses we have that
		$\lf M$ is a~non-zero $\lf$-subbimodule of $M$ and $\lf N$ is a non-zero $\lf$-subbimodule of $N$, and thus we conclude from the irreducibility of either bimodule that $\lf M=M$ and $\lf N=
		N$. Then (a) is an immediate consequence of Corollary~\ref{sym}\,(c), and (b) is an immediate consequence of Corollary~\ref{sym}\,(d).
	\end{proof}

	As promised earlier, we show now that, in general, neither of the truncated tensor products is associative. More precisely, we have the following result:
	\begin{proposition}\label{nonasso}
		For every non-perfect Leibniz algebra $\lf$ there exist
		$\lf$-bimodules $L$, $M$, and $N$ such that
		\[
			(L\ovo M)\ovo N\not\cong L\ovo(M\ovo N)
		      \quad \text{and} \quad
			(L\und M)\und N\not\cong L\und(M\und N)
		\]
		as $\lf$-bimodules.
	\end{proposition}

	\begin{proof}
		Since by hypothesis $\lf$ is not perfect, there exists a~non-zero linear form $\lambda\in\lf^*$ such that $\lambda(\lf\lf)=0$
		(see the argument in Remark~\ref{nonperfect}). Hence, the left
		$\lf$-module $F_{\pm \lambda}$ with underlying vector space
		$\fbb$ and action $x\cdot 1:=\pm\lambda(x)$ is non-trivial. Now choose the $\lf$-bimodules $L:=F_\lambda^s$, $M:=
		F_{-\lambda}^s$, and $N:=F_\lambda^a$. Then we conclude from Theorem~\ref{main} and Corollary~\ref{triv} that
		\[
			(L\ovo M)\ovo N\cong(L\und M)\und N\cong F_0^{s/a}
			\otimes N\cong N\,.
		\]
		On the other hand, we obtain from Theorem~\ref{main} and
		Corollary~\ref{irred}\,(a)
		that
		\[
			L\ovo(M\ovo N)\cong L\und(M\und N)\cong L\und 0=0\,,
		\]
		which completes the proof.
	\end{proof}

	In the last section of this paper we will establish the conclusion of Proposition~\ref{nonasso} for every non-zero Leibniz algebra whose canonical Lie algebra is finite dimensional (see
	Corollary~\ref{nonasso2}).

	We conclude this section by discussing some examples of the case where $M=N=\lf_{\ad}$, where $\lf_{\ad}$ denotes the
	\emph{adjoint $\lf$-bimodule\/} (see~\cite[Example 3.8]{F1}). Note that, in general, the adjoint Leibniz bimodule is neither symmetric nor anti-sym\-metric. In this special case, we have that
	\begin{eqnarray*}
		T(\lf_{\ad},\lf_{\ad})\subseteq T_0 (\lf_{\ad},\lf_{\ad}):=
		(\lf_{\ad})_0\otimes\lf_{\ad}\lf+\lf_{\ad}\lf\otimes(\lf_{\ad})_0\\
		\subseteq\leib(\lf)\otimes\lf\lf+\lf\lf\otimes\leib(\lf)\,,
	\end{eqnarray*}
	and the second inclusion is an equality in case the characteristic of the ground field is not 2 (see~\cite[Example 3.11]{F1}). In all the following examples the right-most term in the above chain of inclusions is properly contained in $\lf_{\ad}\otimes\lf_{\ad}$, and therefore the corresponding truncated tensor products are non-zero.
	\begin{example}
		{Let $\af=\fbb h\oplus\fbb e$ be the solvable left Leibniz algebra over a~field $\fbb$ with multiplication $he=e$ and $hh=
		eh=ee=0$ (see~\cite[Example 2.3]{F1}). Then we have that
		$(\af_{\ad})_0 =\leib(\af)=\fbb e=\af\af$, and therefore}
		\[
			T(\af_{\ad},\af_{\ad})=T_0 (\af_{\ad},\af_{\ad})=\leib(\af)\otimes
			\af\af+\af\af\otimes\leib(\af)=\fbb(e\otimes e)\subsetneqq
			\af_{\ad}\otimes\af_{\ad}\,,
		\]
		{as the vector space on the left-hand side is one-dimensional and the vector space on the right-hand side is four-dimensional.}
	\end{example}

	\begin{example}
		{Let $\nf=\fbb e\oplus\fbb c$ be the nilpotent left and right Leibniz algebra over a~field $\fbb$ with multiplication $ee
		=c$ and $ec=ce=cc=0$ (see~\cite[Example~2.4]{F1}). Then we have that $\leib(\nf)=\fbb c=\nf\nf$, and therefore}
		\[
			\leib(\nf)\otimes\nf\nf+\nf\nf\otimes\leib(\nf)=\fbb
			(c\otimes c)\subsetneqq\nf_{\ad}\otimes\nf_{\ad}\,,
		\]
		{as the dimension of the vector space on the left-hand side is 1 and the dimension of the vector space on the right-hand side is 4.}

		{Because of}
		\begin{eqnarray*}
			(\nf_{\ad})_0 =
			\left\{
			\begin{array}
				{cl}
				\fbb c & \mbox{if }\ch(\fbb)\ne 2\\
				0 & \mbox{if }\ch(\fbb)=2\,,
			\end{array}
			\right.
		\end{eqnarray*}
		{we obtain that}
		\begin{eqnarray*}
			T(\nf_{\ad},\nf_{\ad})=
			\left\{
			\begin{array}
				{cl}
				\fbb(c\otimes c) & \mbox{if }\ch(\fbb)\ne 2\\
				0 & \mbox{if }\ch(\fbb)=2
			\end{array}
			\right.
		\end{eqnarray*}
		and
		\begin{eqnarray*}
			T_0 (\nf_{\ad},\nf_{\ad})=(\nf_{\ad})_0\otimes\nf_{\ad}\nf+\nf_{\ad}\nf
			\otimes(\nf_{\ad})_0 =
			\left\{
			\begin{array}
				{cl}
				\fbb(c\otimes c) & \mbox{if }\ch(\fbb)\ne 2\\
				0 & \mbox{if }\ch(\fbb)=2\,.
			\end{array}
			\right.
		\end{eqnarray*}
		{Hence, we still have that}
		\[
			T(\nf_{\ad},\nf_{\ad})=T_0 (\nf_{\ad},\nf_{\ad})
		\]
		{in any characteristic.}
	\end{example}

	\begin{example}
		{Consider the hemi-semidirect product $\ssf=\slf_2 (\cbb)
		\ltimes_{\hemi} L(1)$, where $\slf_2 (\cbb)$ is the three-dimensional simple Lie algebra of traceless $2\times 2$ matrices with complex coefficients and $L(1)$ is the two-dimensional irreducible $\slf_2
		(\cbb)$-module (see~\cite[Example~2.5]{F1}). Then $\ssf$ is a simple Leibniz algebra (see~\cite[Theorem 2.3]{F2}), and we have that $\leib(\ssf)=L(1)$ and $\ssf\ssf=\ssf$ (for the former see the proof of~\cite[Theorem 2.3]{F2} and for the latter see~\cite[Proposition~7.1]{F1}), and therefore}
		\[
			T_0 (\ssf_{\ad},\ssf_{\ad})=\leib(\ssf)\otimes\ssf\ssf+\ssf\ssf
			\otimes\leib(\ssf)=L(1)\otimes\ssf+\ssf\otimes L(1) \subsetneqq
			\ssf_{\ad}\otimes\ssf_{\ad}\,,
		\]
		{as the dimension of the vector space on the left-hand side is}
		\begin{eqnarray*}
			\dim_\cbb[L(1)\otimes\ssf+\ssf\otimes L(1)] & = &
			2\dim_\cbb[L(1)\otimes\ssf]
			-\dim_\cbb[L(1)\otimes\ssf\cap\ssf\otimes L(1)]\\
			& \le & 2\dim_\cbb[L(1)\otimes\ssf]=20
		\end{eqnarray*}
		{and the dimension of the vector space on the right-hand side is}
		\[
			\dim_\cbb(\ssf\otimes\ssf)=25\,.
		\]
		{By a~straightforward but tedious computation one can show that, in this case, the equality $T(\ssf_{\ad},\ssf_{\ad})=T_0 (\ssf_{\ad},\ssf_{\ad})$ also holds.}
	\end{example}
	It would be very interesting to find general sufficient conditions on $M$ and $N$ which guarantee that $M\ovo N\ne 0$ or
	$M\und N\ne 0$. (Note that $M\und N\ne 0$ implies $M\ovo N
	\ne 0$.)

	\section{Grothendieck rings}\label{grothendieck}

	Let $\ccal$ be an $\fbb$-linear abelian category such that every object has a~composition series of finite length, and let $\K(\ccal)$
	denote the free abelian group generated by the set of isomorphism classes $\irr(\ccal)$ of irreducible (= simple) objects in $\ccal$. Then we assign to every object $X$ in $\ccal$ its class\footnote{Note that two objects in $\ccal$ belong to the same class exactly when they have the same composition factors (up to isomorphism) with the same multiplicities. In particular, isomorphic objects in $\ccal$ belong to the same class, but, of course, unless the objects are irreducible, the converse is not true.} $[X]\in\K(\ccal)$ by
	\[
		[X]
		:=\sum_{S\in\irr(\ccal)}[X:S][S]\,,
	\]
	where $[X:S]$ denotes the (unique) multiplicity of the irreducible object $S$ in a~composition series of $X$ and $[S]\in\K(\ccal)$
	denotes the isomorphism class of $S$ (see~\cite[Definition 1.5.8]{EGNO}). By virtue of the Jordan-H\"older theorem, the above sum is always finite, and $[Z]= [X]+[Y]$ for all short exact sequences $0\to X\to Z\to Y\to 0$ in
	$\ccal$. In particular, we have that
	\[
		[X\oplus Y]
		=[X]+[Y]\,,
	\]
	i.e., the addition in $\K(\ccal)$ corresponds to taking direct sums in
	$\ccal$.

	Now assume that $\ccal$ is a~multiring category over a~field $\fbb$
	(see~\cite[Definition~4.2.3]{EGNO}). Then the tensor product on $\ccal$
	induces a~multiplication on $\K(\ccal)$ via
	\[
		[S_1 ]
		\cdot[S_2 ]:=[S_1\otimes S_2 ]:=\sum_{[S]\in\irr(\ccal)}[S_1
		\otimes S_2:S][S]
	\]
	for all classes $[S_1 ],[S_2 ]\in\irr(\ccal)$. It follows from~\cite[Lemma
	4.5.1]{EGNO} that the multiplication on $\K(\ccal)$ is associative, and then $\K(\ccal)$ is called the \emph{Grothendieck ring\/} of $\ccal$
	(see~\cite[Definition~4.5.2]{EGNO}).

	Next, we specialize the previous paragraph to the category $\ccal:=
	\mod(\gf)$ of finite-dimensional $\gf$-modules for a~Lie algebra
	$\gf$ and to the category $\ccal:=\fdmod_{\weak}^{\bi}(\lf)$ of finite-dimensional weak $\lf$-bimodules for a~Leibniz algebra
	$\lf$. For brevity we set
	\[
		\K(\gf):=\K(\fdmod(\gf))
	\]
	and
	\[
		\K_{\weak}^{\bi}(\lf):=\K(\fdmod_{\weak}^{\bi}(\lf))\,.
	\]
	Although the category $\fdmod^{\bi}(\lf)$ of finite-dimensional $\lf$-bimodules is not a~monoidal category, it is still a~locally finite $\fbb$-linear abelian category, and therefore we can define the Grothendieck group $\K(\fdmod^{\bi}(\lf))$ in the same way as above which we will denote by $\K^{\bi}(\lf)$. Moreover, the truncated tensor product $\ovo$ from Section~\ref{truncated} induces a~multiplication on $\K^{\bi}(\lf)$ via
	\[
		[M]
		\cdot[N]:=[M\ovo N]:=\sum_{[L]\in\irr^{\bi}(\lf)}
		[M\ovo N:L][L]
	\]
	for all classes $[M],[N]\in\irr^{\bi}(\lf)$, where $\irr^{\bi}(\lf)$
	denotes the set of isomorphism classes of finite-dimensional irreducible $\lf$-bimodules. Since irreducible Leibniz bimodules are either symmetric or anti-symmetric, it follows from
	Theorem~\ref{main} that $M\ovo N=M\und N$ for irreducible
	$\lf$-bimodules $M$ and $N$, and therefore instead of $\ovo$
	we could also have used $\und$ in the definition of the multiplication of $\K^{\bi}(\lf)$, i.e., for the Grothendieck ring $\K^{\bi}(\lf)$
	it does not matter which of the truncated tensor products we use.

	As neither of the truncated tensor products is associative
	(see Proposition~\ref{nonasso}), we cannot expect the
	Grothendieck ring $\K^{\bi}(\lf)$ to be associative, and indeed, we will see later in this section that this is quite often not the case (see Proposition~\ref{nonassoK_0}). In fact, from the latter we will deduce Corollary~\ref{nonasso2}
	which complements Proposition~\ref{nonasso}.

	On the other hand, it is an immediate consequence of Proposition~\ref{comm} and Corollary~\ref{triv} that $\K^{\bi}(\lf)$ is commutative and unital:
	\begin{proposition}\label{unitalcomm}
		The Grothendieck ring $\K^{\bi}(\lf)$ of a~Leibniz algebra $\lf$
		is commutative. Moreover, the class $[F_0^{s/a}]$ of the one-dimensional trivial $\lf$-bimodule $F_0^{s/a}$ is the unity of
		$\K^{\bi}(\lf)$.
	\end{proposition}
	We continue by describing the algebraic structure of the Grothendieck ring $\K^{\bi}(\lf)$ of a~Leibniz algebra $\lf$ as a~certain amalgamated product of two copies of the Grothendieck ring $\K(\lf_{\lie})$ of the canonical Lie algebra $\lf_{\lie}$ associated with $\lf$, which we call the unital commutative product.

	Let $A$ and $B$ be unital commutative (but not necessarily associative)
	rings whose underlying additive abelian groups are free, and let $1_A$
	resp.\ $1_B$ denote the unity of the ring $A$ resp.\ $B$. Furthermore, suppose that $A_1$ is a~basis of the free abelian additive group of $A$
	such that $1_A\in A_1$ and $B_1$ is a~basis of the free abelian additive group of $B$ such that $1_B\in B_1$. Then we define the \emph{unital commutative product\/} $A\circledast B$ of $A$ and $B$ as follows:

	Set $\widetilde{A}:=A_1\setminus\{1_A\}$ and $\widetilde{B}:=B_1
	\setminus\{1_B\}$. Adjoin to the set $\widetilde{A}\cup\widetilde{B}$
	a new unity $1$ and consider the free abelian group $G(A,B)$ generated by this set. Then define a~multiplication on the basis of $G(A,B)$ by
	\begin{gather}\label{A}
		a*a':=aa'\mbox{ for all }a,a'\in\widetilde{A} \\
		\label{B}
		b*b':=bb'\mbox{ for all }b, b'\in\widetilde{B} \\
		\label{zero}
		a*b:=b*a:=0\mbox{ for all }a\in\widetilde{A}\,b\in\widetilde{B} \\
		\label{unity}
		1*x:=x*1:=x\mbox{ for every }x\in\widetilde{A}\cup\widetilde{B}\cup\{1\}
	\end{gather}
	and extend it biadditively to all of $G(A,B)$. While doing so, it is understood that in case $aa'=n1_A +\sum_{\tilde{a}\in\widetilde{A}}n_{\tilde{a}}
	\tilde{a}$ in $A$ for some non-zero integer $n$, we require $a*a'=n1+
	\sum_{\tilde{a}\in\widetilde{A}}n_{\tilde{a}}\tilde{a}$ in $A\circledast
	B$, and similarly for $B$\footnote{It is understood that at most finitely many of the integers $n_{\tilde{a}}$ are non-zero.}.

	Observe that $\widetilde{A}\cup\{1\}$ generates a~subring of $A\circledast
	B$ isomorphic to $A$ and that $\widetilde{B}\cup\{1\}$ generates a~subring of $A\circledast B$ isomorphic to $B$.

	It is clear from the definition of the multiplication that $A\circledast B$ is a~unital commutative ring with unity $1$:
	\begin{proposition}
		If $A$ and $B$ are unital commutative rings, then the unital commutative product $A\circledast B$ is a~unital commutative ring with unity $1$.
	\end{proposition}
	Note that the Grothendieck ring $\K^{\bi}(0)$ of the zero algebra is the ring of integers which is also the smallest example of a~unital commutative product $\zbb\circledast\zbb$. In particular, this unital commutative product is associative. In the following we will see that this is an exception to the rule.

	In fact, it turns out that quite often unital commutative products do not satisfy any of the usual associativity properties (associative, alternative, Jordan identity, power-associative). Note that every commutative ring clearly is reflexive and
	3$^\textrm{rd}$ power\--as\-so\-cia\-tive. Moreover, a~commutative ring is left alternative exactly when it is right alternative. Consequently, a~commutative ring that is not alternative is neither left nor right alternative\footnote{See~\cite{BMS} or~\cite{S}
	for the precise definitions of the associativity properties.}. For later applications, we need the following criterion, which gives sufficient conditions for the failure of each of the first three of the associativity properties mentioned above.

	\begin{lemma}\label{criterion}
		Let $A$ and $B$ be unital commutative rings. Then the following statements hold:
		\begin{enumerate}
			\item[\textrm{(a)}] If there exist elements $a,a'\in\widetilde{A}$ such that $aa'
			=n1_A +\sum_{\tilde{a}\in\widetilde{A}}n_{\tilde{a}}\tilde{a}$
			for some non-zero integer $n$ and $\widetilde{B}\ne\emptyset$, then $A\circledast B$ is not associative.
			\item[\textrm{(b)}] If there exists an element $a\in\widetilde{A}$ such that $aa=n1_A
			+\sum_{\tilde{a}\in\widetilde{A}}n_{\tilde{a}}\tilde{a}$ for some non-zero integer $n$ and $\widetilde{B}\ne\emptyset$, then
			$A\circledast B$ is not alternative.
			\item[\textrm{(c)}] If there exist elements $a\in\widetilde{A}$ and $b\in\widetilde{B}$
			such that $aa=m1_A +\sum_{\tilde{a}\in\widetilde{A}}
			m_{\tilde{a}}\tilde{a}$ and at least one of the integers $m_{\tilde{a}}
			\ne 0$ is non-zero as well as $bb=n1_B +\sum_{\tilde{b}\in
			\widetilde{B}}n_{\tilde{b}}\tilde{b}$ for some non-zero integer $n$, then $A\circledast B$ is not a~Jordan ring.
		\end{enumerate}
	\end{lemma}

	\begin{proof}
		(a): Since by hypothesis $\widetilde{B}\ne\emptyset$, there exists an element
		$b\in\widetilde{B}$. Then one can build the following expressions:
		\[
			(a*a')*b=\left(n1+\sum_{\tilde{a}\in\widetilde{A}}n_{\tilde{a}}\tilde{a}
			\right)*b=nb\ne 0
		\]
		and
		\[
			a*(a'*b)=a*0=0\,,
		\]
		which contradicts associativity.

		The proof of (b) is very similar to the proof of (a) and is left to the reader.

		(c): If we set $a^2:=aa$, we obtain, similar to the proof of (a), that
		\[
			(a^2 *b)*b=(mb)*b=mbb=(mn)1_B +\sum_{\tilde{b}\in\widetilde{B}}
			(mn_{\tilde{b}})\tilde{b}\,.
		\]
		On the other hand, we have that
		\begin{eqnarray*}
			a^2 *(b*b) & = & a^2 *(bb)=\left(m1_A +\sum_{\tilde{a}\in\widetilde{A}}
			m_{\tilde{a}}\tilde{a}\right)*\left(n1_B +\sum_{\tilde{b}\in\widetilde{B}}
			n_{\tilde{b}}\tilde{b}\right)\\
			& = & (mn)1+\sum_{\tilde{b}\in\widetilde{B}}(mn_{\tilde{b}})\tilde{b}
			+\sum_{\tilde{a}\in\widetilde{A}}(nm_{\tilde{a}})\tilde{a}\,.
		\end{eqnarray*}
		Suppose now that $A\circledast B$ is a~Jordan ring. Then $(a^2 *b)*b=
		a^2 *(b*b)$ implies that $nm_{\tilde{a}}=0$ for every $\tilde{a}\in
		\widetilde{A}$, which contradicts the hypothesis.
	\end{proof}

	Recall that every irreducible $\lf$-bimodule is either the symmetrization $M^s$
	or the anti-symmetrization $M^a$ of some irreducible $\lf_{\lie}$-module $M$. Hence,
	\[
		\{[M^s ]\mid[M]\in\irr(\lf_{\lie})\}\cup\{[M^a ]\mid[M]\in\irr(\lf_{\lie})\}
	\]
	is a~basis of the additive abelian group underlying the ring $\K^{\bi}(\lf)$, where $\irr(\lf_{\lie})$ denotes the set of isomorphism classes of finite-dimensional irreducible $\lf_{\lie}$-modules.
	\begin{theorem}\label{communitalprod}
		The Grothendieck ring $\K^{\bi}(\lf)$ of a~Leibniz algebra $\lf$ is isomorphic to the unital commutative product of two copies of the Grothendieck ring $\K
		(\lf_{\lie})$ of the canonical Lie algebra $\lf_{\lie}$ associated with $\lf$, i.e.,
		\[
			\K^{\bi}(\lf)\cong\K(\lf_{\lie})\circledast\K(\lf_{\lie})
		\]
		as unital commutative rings.
	\end{theorem}

	\begin{proof}
		Let $A$ denote the free abelian subgroup of $\K^{\bi}(\lf)$ that generated by
		\[\{[M^s ]\mid[M]\in\irr(\lf_{\lie})\},\] and let $B$ denote the free abelian subgroup of $\K^{\bi}(\lf)$ that is generated by \[\{[M^a ]\mid[M]\in\irr(\lf_{\lie})\}.\]

		It follows from Corollary~\ref{sym}\,(a) and (b) that
		\[
			M^s\ovo N^s =M^s\otimes N^s =(M\otimes N)^s
		\]
		and
		\[
			M^a\ovo N^a =M^a\otimes N^a =(M\otimes N)^a
		\]
		for all irreducible $\lf_{\lie}$-modules $M$ and $N$. This shows that $\K^{\bi}(\lf)$
		contains $A$ and $B$ as unital subrings. In particular, the relations
		\eqref{A} and \eqref{B} hold in $\K^{\bi}(\lf)$. Moreover, we have that $A\cong\K
		(\lf_{\lie})$ via the isomorphism $[M]\mapsto[M^s ]$ and $B\cong\K(\lf_{\lie})$ via the isomorphism $[M]\mapsto[M^a ]$.

		In addition, we obtain from Corollary~\ref{irred} that
		\[
			M^s\ovo N^a =M^a\ovo N^s =0
		\]
		for all non-trivial irreducible $\lf_{\lie}$-modules $M$ and $N$. This implies that
		\eqref{zero} holds in $\K^{\bi}(\lf)$.

		Finally, \eqref{unity} is an immediate consequence of Corollary~\ref{triv}.
	\end{proof}

	In particular, Theorem~\ref{communitalprod} shows that the Grothendieck ring
	$\K^{\bi}(\lf)$ only depends on the canonical Lie algebra $\lf_{\lie}$. By applying
	Theorem~\ref{communitalprod} both to $\lf$ and $\lf_{\lie}$, we obtain the following result, which shows that for computing the Grothendieck ring $\K^{\bi}
	(\lf)$ it is enough to consider only Lie algebras:
	\begin{corollary}\label{reduc}
		The Grothendieck ring $\K^{\bi}(\lf)$ of a~Leibniz algebra $\lf$ is isomorphic to the Grothendieck ring $\K^{\bi}(\lf_{\lie})$ of the canonical Lie algebra
		$\lf_{\lie}$ associated with $\lf$, i.e.,
		\[
			\K^{\bi}(\lf)\cong\K^{\bi}(\lf_{\lie})
		\]
		as unital commutative rings.
	\end{corollary}
	Next, we illustrate Theorem~\ref{communitalprod} and Corollary~\ref{reduc}
	by the following example.
	\begin{example}\label{1dim}
		{Let $\ef:=\fbb e$ be the one-dimensional Lie algebra over an algebraically closed field $\fbb$ of arbitrary characteristic. Then the finite-dimensional irreducible $\ef$-modules are the one-dimensional modules $F_\lambda$, where $e$ acts as multiplication by a~scalar
		$\lambda$ on the ground field $\fbb$, and therefore the isomorphism classes of the finite-dimensional irreducible $\ef$-bimodules are}
		\[
			\{[F_\lambda^s ]\mid\lambda\in\fbb\}\cup\{[F_\mu^a ]\mid\mu\in\fbb\}\,.
		\]
		{By virtue of Corollary~\ref{sym}\,(a), (b) and Corollary~\ref{irred}, the truncated tensor products are either the ``natural'' tensor product defined in Section~\ref{weak} or zero, i.e.,}
		\[
			F_\lambda^s\ovo F_\mu^s\cong F_{\lambda+\mu}^s\,,\quad\quad
			F_\lambda^a\ovo F_\mu^a\cong F_{\lambda+\mu}^a\,,\quad\quad
			F_\lambda^s\ovo F_\mu^a\cong F_\mu^a\ovo F_\lambda^s =0\,,
		\]
		{where in the last isomorphism we suppose that $\lambda\not=0$
		and $\mu\not=0$. Moreover, it follows from Corollary~\ref{triv} that the one-dimensional trivial $\ef$-bimodule $F_0^{s/a}:=F_0^s =F_0^a$ acts as an identity for both truncated tensor products. As a~consequence, we obtain that the Grothendieck ring $\K(\ef)=\bigoplus_{\lambda\in
		\fbb}\zbb[F_\lambda]$ of the Lie algebra $\ef$ is isomorphic to the integral group ring $\zbb[\fbb^+ ]$ of the additive group $\fbb^+$ of the ground field $\fbb$, and thus it follows from Theorem~\ref{communitalprod} that the Grothendieck ring $\K^{\bi}(\ef)$ is the unital commutative product of two copies of $\zbb[\fbb^+ ]$:}
		\[
			\K^{\bi}(\ef)\cong\zbb[\fbb^+ ]\circledast\zbb[\fbb^+ ]\,,
		\]
		{where one copy of $\zbb[\fbb^+ ]$ corresponds to the symmetrizations of the finite-dimen\-sional irreducible $\ef$-modules and the other copy corresponds to the anti-symme\-trizations of the finite-dimensional irreducible
		$\ef$-modules.}

		{More generally, it then follows from Corollary~\ref{reduc} that $\K^{\bi}
		(\lf)\cong\zbb[\fbb^+ ]\circledast\zbb[\fbb^+ ]$ for every Leibniz algebra
		$\lf$ whose canonical Lie algebra $\lf_{\lie}$ is one-dimensional. In particular, this applies to every two-dimensional non-Lie Leibniz algebra.}

		{Note that the classes $\pm[F_\lambda^s ]$, $\pm[F_\lambda^a ]$
		($\lambda\in\fbb$) are at the same time multiplicatively invertible and zero divisors. This already shows that $\K^{\bi}(\ef)$ is not associative because in a~unital commutative associative ring multiplicatively invertible elements are never zero divisors and vice versa. But this can also be seen more explicitly as follows (cf.\ also the proof of Proposition~\ref{nonasso}). Indeed, by fixing a non-zero scalar $\lambda\in\fbb$, we have that}
		\[
			([F_\lambda^s ]\cdot[F_{-\lambda}^s ])\cdot[F_\lambda^a ]=
			[F_0^{s/a}]\cdot[F_\lambda^a ]=[F_\lambda^a ]\,.
		\]
		{On the other hand, we have that}
		\[
			[F_\lambda^s ]
			\cdot([F_{-\lambda}^s ]\cdot[F_\lambda^a ])=
			[F_\lambda^s ]\cdot0= 0\,,
		\]
		{which again shows that the multiplication of $\K^{\bi}(\ef)$ is not associative.}

		{However, we will see in Theorem~\ref{alt} and Corollary~\ref{alt2} below
		that $\K^{\bi}(\ef)$ is an alternative power-associative Jordan ring.}
	\end{example}
	By employing Lie's theorem (see~\cite[Theorem 4.1]{H}), we can generalize
	Example~\ref{1dim} considerably for ground fields of characteristic zero.
	\begin{theorem}\label{lie}
		Let $\gf$ be a~finite-dimensional solvable Lie algebra over an algebraically closed field $\fbb$ of characteristic zero. Then
		\[
			\K(\gf)\cong\zbb[\fbb^+\times\cdots\times\fbb^+ ]
		\]
		as unital commutative rings, where $\fbb^+$ denotes the additive group of $\fbb$ and the number of factors $\fbb^+$ in the group ring is
		$\dim_\fbb\gf/[\gf,\gf]$\footnote{Here $[\gf,\gf]:=\langle [x,y]\in\gf\mid x,y\in\gf\rangle_\fbb$ is the \emph{derived subalgebra\/} of the Lie algebra $\gf$.}.
	\end{theorem}

	\begin{proof}
		According to Lie's theorem, every finite-dimensional irreducible
		$\gf$-module is one-dimensional. Since the derived subalgebra
		$[\gf,\gf]$ of $\gf$ acts trivially on every one-dimensional $\gf$-module, the finite-dimensional irreducible $\gf$-modules are in bijection with
		$(\gf/[\gf,\gf])^*$ from which the assertion follows similarly to the argument in Example~\ref{1dim}.
	\end{proof}

	As an immediate consequence of Theorem~\ref{communitalprod} and Theorem~\ref{lie} we obtain the following result:
	\begin{corollary}\label{solv}
		Let $\lf$ be a~solvable Leibniz algebra over an algebraically closed field $\fbb$
		of characteristic zero whose canonical Lie algebra $\lf_{\lie}$ is finite dimensional. Then
		\[
			\K^{\bi}(\lf)\cong\zbb[\fbb^+\times\cdots\times\fbb^+ ]\circledast\zbb
			[\fbb^+\times\cdots\times\fbb^+ ]
		\]
		as unital commutative rings, where $\fbb^+$ denotes the additive group of $\fbb$ and the number of factors $\fbb^+$ in each of the group rings is
		$\dim_\fbb\lf_{\lie}/[\lf_{\lie},\lf_{\lie}]$.
	\end{corollary}
	Next, we prove that the Grothendieck ring $\K^{\bi}(\lf)$ of a~Leibniz algebra $\lf$
	quite often is not associative. For this we need the following result that should be well-known but for which we could not find a~reference.
	\begin{lemma}\label{nontrivirred}
		Every non-zero finite-dimensional Lie algebra has a~finite-dimensional non-trivial irreducible module.
	\end{lemma}

	\begin{proof}
		Let $\gf$ be a~non-zero finite-dimensional Lie algebra and suppose that the trivial module is the only finite-dimensional irreducible $\gf$-module. But then it follows from~\cite[Proposition 1]{F0} that $\gf$ is nilpotent. As $\gf\ne 0$,
		$\gf$ is not perfect, and thus there exists a~non-zero linear form $\lambda
		\in\gf^*$ such that $\lambda([\gf,\gf])=0$ (see the argument in Remark~\ref{nonperfect}). Hence, the $\gf$-module $F_\lambda$ with underlying vector space $\fbb$ and $\gf$-action $x\cdot 1:=\lambda(x)$ clearly defines a non-trivial finite-dimensional irreducible $\gf$-module contradicting our assumption.
	\end{proof}

	\begin{proposition}\label{nonassoK_0}
		Let $\lf$ be a~non-zero Leibniz algebra whose canonical Lie algebra
		$\lf_{\lie}$ is finite dimensional. Then the Grothendieck ring $\K^{\bi}(\lf)$ is not associative.
	\end{proposition}

	\begin{proof}
		Since $\lf\ne 0$, we obtain from~\cite[Proposition 2.20]{F1}
		that $\lf_{\lie}\ne 0$. Then it follows from Lemma~\ref{nontrivirred} that there exists a~finite-dimensional non-trivial irreducible $\lf_{\lie}$-module $M$, and thus $M^*$ is also a~non-trivial irreducible $\lf_{\lie}$-module. Hence, the contraction $M^*\otimes M\to F_0$, $\mu\otimes m\mapsto\mu(m)$ is an epimorphism of $\lf_{\lie}$-modules, where $F_0$ denotes the one-dimensional trivial $\lf_{\lie}$-module. In particular, we have that $[M^*\otimes M:F_0 ]\ne
		0$, and we deduce from the definition of the multiplication in $\K(\lf_{\lie})$
		that
		\[
			[M^* ]
			\cdot[M]=\sum_{[L]\in\irr(\lf_{\lie})}[M^*\otimes M:L][L]\,.
		\]
		Hence, the assertion follows from Theorem~\ref{communitalprod}
		and Lemma~\ref{criterion}\,(a).
	\end{proof}

	As an immediate consequence of Proposition~\ref{nonassoK_0} we obtain the next result which complements Proposition~\ref{nonasso} (see Example~\ref{nonasso3} below):
	\begin{corollary}\label{nonasso2}
		Let $\lf$ be a~non-zero Leibniz algebra whose canonical Lie algebra
		$\lf_{\lie}$ is finite dimensional. Then there exist
		$\lf$-bimodules $L$, $M$, and $N$ such that
		\[
			(L\ovo M)\ovo N\not\cong L\ovo(M\ovo N)
		\]
		and
		\[
			(L\und M)\und N\not\cong L\und(M\und N)
		\]
		as $\lf$-bimodules.
	\end{corollary}

	\begin{example}\label{nonasso3}
		{Let $\slf_2 (\cbb)$ be the three-dimensional simple Lie algebra of traceless complex $2\times 2$ matrices, and let $M(\lambda)$ denote the Verma module of highest weight $\lambda$ which is an infinite-dimensional
		$\slf_2 (\cbb)$-module. (Here we identify every complex multiple of the unique fundamental weight with its coefficient.) It is well known (see
		\cite[Exercise 7\,(c) in
		Section 7]{H}) that $M(\lambda)$ is irreducible if $\lambda+1$
		is not a~dominant integral weight (i.e., with our identification, $\lambda+1$
		is not a~non-negative integer). Hence, it follows from~\cite[Theorem 2.3]{F2}
		in conjunction with~\cite[Proposition 7.1]{F1} that the hemi-semidirect product $\slf_2 (\cbb)\ltimes_{\hemi} M(-2)$ is a~perfect Leibniz algebra whose canonical Lie algebra is finite dimensional.\footnote{This example also shows that the hypothesis in Corollary~\ref{nonasso2} (and in several other results in this section) is weaker than just to assume that the Leibniz algebra itself is finite dimensional.} Consequently, Corollary~\ref{nonasso2}
		applies, but Proposition~\ref{nonasso} does not.}

		{On the other hand, for every infinite-dimensional non-perfect Lie algebra (for example, the infinite-dimensional Heisenberg algebra)
		Proposition~\ref{nonasso} applies, but Corollary~\ref{nonasso2} does not.}
	\end{example}
	Similar to the proof of Corollary~\ref{solv}, we can employ Lie's theorem to show that the Grothendieck ring $\K^{\bi}(\lf)$ of a~finite-dimensional solvable
	Leibniz algebra $\lf$ over an algebraically closed field of characteristic zero is alternative, or even slightly more general:
	\begin{theorem}\label{alt}
		Let $\lf$ be a~solvable Leibniz algebra over an algebraically closed field of characteristic zero whose canonical Lie algebra $\lf_{\lie}$ is finite dimensional. Then the Grothendieck ring $\K^{\bi}(\lf)$ is alternative.
	\end{theorem}

	\begin{proof}
		Since by Proposition~\ref{unitalcomm}, $\K^{\bi}(\lf)$ is commutative, it suffices to prove that $\K^{\bi}(\lf)$ is left alternative. By virtue of Lie's theorem, every finite-dimensional irreducible $\lf_{\lie}$-module is one-dimensional. Hence, we have that
		\[
			\irr(\lf_{\lie})=\{[F_\lambda]\mid\lambda\in\Lambda\}\,,
		\]
		where
		\[
			\Lambda:=\{\lambda\in\lf_{\lie}^*\mid\lambda([\lf_{\lie},\lf_{\lie}])=0\}\,,
		\]
		and therefore
		\[
			\irr^{\bi}(\lf)=\{[F_0^{s/a}]\}\cup\{[F_\lambda^s ]\mid\lambda\in\Lambda\setminus\{0\}\}
			\cup\{[F_\mu^a ]\mid\lambda\in\Lambda\setminus\{0\}\}\,.
		\]
		Let
		\[
			u=m_0 [F_0^{s/a}]+\sum_{\lambda\in\Lambda\setminus\{0\}}m_\lambda[F_\lambda^s ]
			+\sum_{\mu\in\Lambda\setminus\{0\}}m_\mu[F_\mu^a ]
		\]
		and
		\[
			v=n_0 [F_0^{s/a}]+\sum_{\nu\in\Lambda\setminus\{0\}}n_\nu[F_\nu^s ]
			+\sum_{\eta\in\Lambda\setminus\{0\}}n_\eta[F_\eta^a ]
		\]
		be arbitrary elements in $\K^{\bi}(\lf)$. Then we have that
		\begin{eqnarray*}
			u^2 & = & m_0^2 [F_0^{s/a}]+\sum_{\lambda\in\Lambda\setminus\{0\}}(2m_0 m_\lambda)
			[F_\lambda^s ]+\sum_{\mu\in\Lambda\setminus\{0\}}(2m_0 m_\mu)[F_\mu^a ]\\
			&& +\sum_{\lambda',\lambda''\in\Lambda\setminus\{0\}}(m_{\lambda'}m_{\lambda''})
			[F_{\lambda'+\lambda''}^s ]+\sum_{\mu',\mu''\in\Lambda\setminus\{0\}}
			(m_{\mu'}m_{\mu''})[F_{\mu'+\mu''}^a ]\,,
		\end{eqnarray*}
		and thus
		\begin{eqnarray*}
			u^2\cdot v & = & (m_0^2 n_0 )[F_0^{s/a}]+\sum_{\nu\in\Lambda\setminus\{0\}}(m_0^2 n_\nu)
			[F_\nu^s ]+\sum_{\eta\in\Lambda\setminus\{0\}}(m_0^2 n_\eta)[F_\eta^a ]\\
			&& +\sum_{\lambda\in\Lambda\setminus\{0\}}(2m_0 m_\lambda n_0 )[F_\lambda^s ]
			+\sum_{\lambda,\nu\in\Lambda\setminus\{0\}}(2m_0 m_\lambda n_\nu)[F_{\lambda+\nu}^s ]\\
			&& +\sum_{\mu\in\Lambda\setminus\{0\}}(2m_0 m_\mu n_0 )[F_\mu^a ]
			+\sum_{\mu,\eta\in\Lambda\setminus\{0\}}(2m_0 m_\mu n_\eta)[F_{\mu+\eta}^a ]\\
			&& +\sum_{\lambda',\lambda''\in\Lambda\setminus\{0\}}(m_{\lambda'}m_{\lambda''}n_0 )
			[F_{\lambda'+\lambda''}^s ]+\sum_{\lambda',\lambda'',\eta\in\Lambda\setminus\{0\}}
			(m_{\lambda'}m_{\lambda''}n_\nu)[F_{\lambda'+\lambda''+\nu}^s ]\\
			&& +\sum_{\mu',\mu''\in\Lambda\setminus\{0\}}(m_{\mu'}m_{\mu''}n_0 )[F_{\mu'+\mu''}^a ])
			+\sum_{\mu',\mu'',\eta\in\Lambda\setminus\{0\}}(m_{\mu'}m_{\mu''}n_\eta)[F_{\mu'+\mu''+\eta}^a ])\,.
		\end{eqnarray*}
		On the other hand, we have that
		\begin{eqnarray*}
			u\cdot v & = & (m_0 n_0 )[F_0^{s/a}]+\sum_{\nu\in\Lambda\setminus\{0\}}(m_0 n_\nu)[F_\nu^s ]
			+\sum_{\eta\in\Lambda\setminus\{0\}}(m_0 n_\eta)[F_\eta^a ]\\
			&& +\sum_{\lambda\in\Lambda\setminus\{0\}}(m_\lambda n_0 )[F_\lambda^s ]
			+\sum_{\lambda,\nu\in\Lambda\setminus\{0\}}(m_\lambda n_\nu)[F_{\lambda+\nu}^s ]\\
			&& +\sum_{\mu\in\Lambda\setminus\{0\}}(m_\mu n_0 )[F_\mu^a ]
			+\sum_{\mu,\eta\in\Lambda\setminus\{0\}}(m_\mu n_\eta)[F_{\mu+\eta}^a ]\,,
		\end{eqnarray*}
		and then
		\begin{eqnarray*}
			u\cdot(u\cdot v) & = & (m_0^2 n_0 )[F_0^{s/a}]+\sum_{\nu\in\Lambda\setminus\{0\}}(m_0^2 n_\nu)
			[F_\nu^s ]+\sum_{\eta\in\Lambda\setminus\{0\}}(m_0^2 n_\eta)[F_\eta^a ]\\
			&& +\underbrace{\sum_{\lambda\in\Lambda\setminus\{0\}}(m_0 m_\lambda n_0 )[F_\lambda^s ]}_{1}
			+\underbrace{\sum_{\lambda,\nu\in\Lambda\setminus\{0\}}(m_0 m_\lambda n_\nu)[F_{\lambda+\nu}^s ]}_{2}\\
			&& +\underbrace{\sum_{\mu\in\Lambda\setminus\{0\}}(m_0 m_\mu n_0 )[F_\mu^a ]}_{3}
			+\underbrace{\sum_{\mu,\eta\in\Lambda\setminus\{0\}}(m_0 m_\mu n_\eta)[F_{\mu+\eta}^a ]}_{4}\\
			&& +\underbrace{\sum_{\lambda\in\Lambda\setminus\{0\}}(m_\lambda m_0 n_0 )[F_\lambda^s ]}_{1}
			+\underbrace{\sum_{\lambda,\nu\in\Lambda\setminus\{0\}}(m_\lambda m_0 n_\nu)[F_{\lambda+\nu}^s ]}_{2}\\
			&& +\sum_{\lambda',\lambda''\in\Lambda\setminus\{0\}}(m_{\lambda'}m_{\lambda''}n_0 )
			[F_{\lambda'+\lambda''}^s ]
			+\sum_{\lambda',\lambda'',\nu\in\Lambda\setminus\{0\}}(m_{\lambda'}m_{\lambda''}n_\nu)
			[F_{\lambda'+\lambda''+\nu}^s ]\\
			&& +\underbrace{\sum_{\mu\in\Lambda\setminus\{0\}}(m_\mu m_0 n_0 )[F_\mu^a ]}_{3}
			+\underbrace{\sum_{\mu,\eta\in\Lambda\setminus\{0\}}(m_\mu m_0 n_\eta)[F_{\mu+\eta}^a ]}_{4}\\
			&& +\sum_{\mu',\mu''\in\Lambda\setminus\{0\}}(m_{\mu'}m_{\mu''}n_0 )[F_{\mu'+\mu''}^a ]
			+\sum_{\mu', \mu'', \eta\in\Lambda\setminus\{0\}}(m_{\mu'}m_{\mu''}n_\eta)[F_{\mu'+\mu''+\eta}^a ]\\
			& = & (m_0^2 n_0 )[F_0^{s/a}]+\sum_{\nu\in\Lambda\setminus\{0\}}(m_0^2 n_\nu)[F_\nu^s ]
			+\sum_{\eta\in\Lambda\setminus\{0\}}(m_0^2 n_\eta)[F_\eta^a ]\\
			&& +\sum_{\lambda\in\Lambda\setminus\{0\}}(2m_0 m_\lambda n_0 )[F_\lambda^s ]
			+\sum_{\lambda,\nu\in\Lambda\setminus\{0\}}(2m_0 m_\lambda n_\nu)[F_{\lambda+\nu}^s ]\\
			&& +\sum_{\mu\in\Lambda\setminus\{0\}}(2m_0 m_\mu n_0 )[F_\mu^a ]
			+\sum_{\mu,\eta\in\Lambda\setminus\{0\}}(2m_0 m_\mu n_\eta)[F_{\mu+\eta}^a ]\\
			&& +\sum_{\lambda',\lambda''\in\Lambda\setminus\{0\}}(m_{\lambda'}m_{\lambda''}n_0 )
			[F_{\lambda'+\lambda''}^s ]
			+\sum_{\lambda',\lambda'',\nu\in\Lambda\setminus\{0\}}(m_{\lambda'}m_{\lambda''}n_\nu)
			[F_{\lambda'+\lambda''+\nu}^s ]\\
			&& +\sum_{\mu',\mu''\in\Lambda\setminus\{0\}}(m_{\mu'}m_{\mu''}n_0 )[F_{\mu'+\mu''}^a ]
			+\sum_{\mu', \mu'', \eta\in\Lambda\setminus\{0\}}(m_{\mu'}m_{\mu''}n_\eta)[F_{\mu'+\mu''+\eta}^a ]\,,
		\end{eqnarray*}
		which completes the proof.
	\end{proof}

	The next example shows that Theorem~\ref{alt} is not true in non-zero characteristic:
	\begin{example}\label{primchar}
		{Consider the non-abelian solvable restricted Lie algebra $L=\fbb t\ltimes I$ in
		Theorem 3.1 of~\cite{FS} over an algebraically closed field $\fbb$ of prime characteristic
		$p$ and choose $Z$ to be properly contained in $T$. Then it follows from Corollary~1.2\,b)
		in~\cite{F-1} that each of the $p$-dimensional irreducible restricted modules is self-dual. Hence, we can argue similar to the proof of Proposition~\ref{nonassoK_0}
		by applying part (b) instead of part (a) of Lemma~\ref{criterion} to show that the
		Grothendieck ring $\K^{\bi}(L)$ is not alternative.}
	\end{example}
	For the convenience of the reader we include the following well-known result:
	\begin{lemma}\label{jordanpowerasso}
		Let $R$ be a~commutative ring. Then the following statements hold:
		\begin{enumerate}
			\item[\textrm{(a)}] If $R$ is alternative, then $R$ is a~Jordan ring.
			\item[\textrm{(b)}] If $R$ is a~Jordan ring of characteristic $\ne 2,3,5$, then $R$ is power-associative.
		\end{enumerate}
	\end{lemma}

	\begin{proof}
		(a): Since by hypothesis $R$ is alternative and commutative, we obtain from the left alternative law, i.e., $r^2 s=r(rs)$ holds for arbitrary elements $r,s\in R$, in conjunction with the commutativity of $R$ that
		\[
			(u^2 v)u=[u(uv)]u=u[u(vu)]=u^2 (vu)
		\]
		for all elements $u,v\in R$.

		(b): According to~\cite[Lemmas 3 and 4]{A}, it is enough to prove that $R$ is {\em
		$4^\textrm{rd}$ power-associative\/}, i.e., $R$ satisfies $r^2 r^2 =(r^2 r)r$ for every element $r\in R$. But the latter is a~special case of the Jordan identity.
	\end{proof}

	\begin{remark}
		{Note that Lemma~\ref{jordanpowerasso} suffices for our purposes, but the implications are known to be true more generally. Namely, it follows from a theorem of Emil Artin that every alternative ring is power-associative. Moreover, every (not necessarily commutative) Jordan ring is power-associative (see Fact
		6 on p.\ 19 in~\cite{BMS}).}
	\end{remark}
	As an immediate consequence of Theorem~\ref{alt} in conjunction with
	Proposition~\ref{unitalcomm} and Lemma~\ref{jordanpowerasso} we obtain the following result:
	\begin{corollary}\label{alt2}
		Let $\lf$ be a~solvable Leibniz algebra over an algebraically closed field of characteristic zero whose canonical Lie algebra $\lf_{\lie}$ is finite dimensional. Then the Grothendieck ring $\K^{\bi}(\lf)$ is a~power-associative commutative Jordan ring.
	\end{corollary}
	Now, we consider the Grothendieck rings $\K^{\bi}(\lf)$ of finite-dimensional semi-simple Leibniz algebras $\lf$ and begin with the smallest possible example:
	\begin{example}
		{Let $\slf_2 (\cbb)$ denote the three-dimensional simple Lie algebra of traceless complex $2\times 2$ matrices. Note that the category
		$\fdmod^{\bi} (\slf_2 (\cbb))$ of finite-dimensional $\slf_2 (\cbb)$-bimodules has already been studied in~\cite{LP2}. Here we consider the Grothendieck ring
		$\K^{\bi}(\slf_2 (\cbb))$ of this category.}

		{It is well known that the isomorphism classes of the finite-dimensional irreducible $\slf_2 (\cbb)$-modules are given by}
		\[
			\{[L(n)]\mid\dim_\cbb L(n)=n+1\,,\,n\in\nbb_0\}
		\]
		{(see~\cite[Theorem 7.2]{H}), and therefore the isomorphism classes of finite-dimensional irreducible $\slf_2 (\cbb)$-bimodules are}
		\[
			\{[L(n)^s ]\mid n\in\nbb_0\}\cup\{[L(n)^a ]\mid n\in\nbb_0\}\,.
		\]
		{By virtue of Corollary~\ref{sym}\,(a), (b) and Corollary~\ref{irred}, the truncated tensor products are either the ``natural'' tensor product defined in
		Section~\ref{weak} or zero, i.e.,}
		\begin{gather*}
			L(m)^s\ovo L(n)^s\cong[L(m)\otimes L(n)]^s,\quad\quad L(m)^a\ovo L(n)^a
			\cong[L(m)\otimes L(n)]^a\,, \\
			L(m)^s\ovo L(n)^a\cong L(n)^a\ovo L(m)^s =0\,,
		\end{gather*}
		{where in the last isomorphism we assume that $m>0$ and $n>0$. Moreover, the tensor products $L(m)\otimes L(n)$ can be obtained from the Clebsch-Gordan formula (see~\cite[Exercise~7 in Section 22]{H}). Finally, it follows from Corollary~\ref{triv} that the one-dimensional trivial $\slf_2 (\cbb)$-bimodule $L(0)^{s/a}:=
		L(0)^s =L(0)^a$ acts as an identity for both truncated tensor products.}

		{As a~consequence of the Clebsch-Gordan formula, we obtain that the
		Grothen\-dieck ring $\K(\slf_2 (\cbb))$ of the Lie algebra $\slf_2 (\cbb)$ is isomorphic to the polynomial ring $\zbb[t]$ in one variable generated by the class $[L(1)]$
		of the two-dimensional irreducible $\slf_2 (\cbb)$-module, and therefore it follows from
		Theorem~\ref{communitalprod} that the Grothendieck ring $\K^{\bi}(\slf_2 (\cbb))$
		is the unital commutative product of two copies of $\zbb[t]$:}
		\[
			\K^{\bi}(\slf_2 (\cbb))\cong\zbb[t]\circledast\zbb[t]\,,
		\]
		{where one copy of $\zbb[t]$ corresponds to the symmetrizations of the finite-dimensional irreducible $\slf_2 (\cbb)$-modules and the other copy corresponds to the anti-symmetri\-zations of the finite-dimensional irreducible
		$\slf_2 (\cbb)$-modules.}

		{More generally, it then follows from Corollary~\ref{reduc} that $\K^{\bi}(\lf)
		\cong\zbb[t]\circledast\zbb[t]$ for every Leibniz algebra $\lf$ whose canonical
		Lie algebra $\lf_{\lie}$ is isomorphic to $\slf_2 (\cbb)$. In particular, this applies to the simple\footnote{See Theorem 2.3 in~\cite{F2}.} non-Lie Leibniz algebras
		$\slf_2 (\cbb)\ltimes_{\hemi} L(n)$ for every positive integer $n$.}

		{Note also that $\K^{\bi}(\slf_2 (\cbb))$ is neither alternative nor a~Jordan ring. Indeed, we obtain from the Clebsch-Gordan formula that}
		\[
			L(m)\otimes L(n)\cong L(m+n)\oplus L(m+n-2)\oplus\dotsb\oplus L(m-n)
		\]
		{for all integers $m\geq n\ge 0$, and where exactly $n+1$ summands occur on the right-hand side. In particular, $L(1)\otimes L(1)\cong L(0)\oplus
		L(2)$ contains the one-dimensional trivial module and a~non-trivial irreducible module as direct summands. In view of Lemma~\ref{criterion}\,(c), this implies in turn that $\K^{\bi}(\slf_2 (\cbb))$ is not a~Jordan ring, and therefore it follows from Lemma~\ref{jordanpowerasso}\,(a) that $\K^{\bi}(\slf_2 (\cbb))$ is also not alternative.}
	\end{example}
	Note that, for dimension reasons, every finite-dimensional irreducible
	$\slf_2 (\cbb)$-mo\-dule is self-dual. In fact, every finite-dimensional non-zero semi-simple Lie algebra over a~field of characteristic zero has a finite-dimensional non-trivial self-dual irreducible module, namely, the adjoint module of one of its minimal ideals, and so we can apply part (b)
	of Lemma~\ref{criterion}. It turns out that even more is true which enables us to apply part (c) of this lemma. We will use this observation to show that the Grothendieck ring $\K^{\bi}(\lf)$ of every finite-dimensional non-zero semi-simple Leibniz algebra over a~field of characteristic zero is neither alternative nor a~Jordan ring which is contrary to the behavior of solvable
	Leibniz algebras as we have seen in Theorem~\ref{alt} and Corollary~\ref{alt2}.
	\begin{theorem}\label{nonalt}
		Let $\lf$ be a~non-zero semi-simple Leibniz algebra over a~field
		$\fbb$ of characteristic zero whose canonical Lie algebra $\lf_{\lie}$
		is finite dimensional. Then the Grothendieck ring $\K^{\bi}(\lf)$ is neither alternative nor a~Jordan ring.
	\end{theorem}

	\begin{proof}
		As in the proof of Proposition~\ref{nonassoK_0}, we have that $\gf:=
		\lf_{\lie}\ne 0$. Note that it follows from~\cite[Proposition 7.8]{F1} that
		$\gf$ is necessarily semi-simple. Now let $M$ be a~minimal ideal of
		$\gf$. Then $M$ is a~simple Lie algebra (see~\cite[Theorem 5.2]{H})
		and a~finite-dimensional non-trivial irreducible $\gf$-module. Since the Killing form of $M$ is $\gf$-invariant (see~\cite[Lemma~5.1]{H})
		and non-degenerate (see~\cite[Theorem~5.1]{H}), $M$ and its linear dual $M^*$ are isomorphic $\gf$-modules (cf.\ Proposition 4.6 in
		Chapter 3 of~\cite{SF}).

		Furthermore, it follows from Weyl's theorem on complete reducibility
		(see~\cite[Theorem~6.3]{H}) and the bijection between finite-dimensional irreducible $\gf$-modules and the set of dominant integral weights
		$\Lambda^+$ with respect to a~Cartan subalgebra of $\gf$ (see~\cite[Corollary~21.2]{H}) that
		\[
			M\otimes M=\bigoplus_{\lambda\in\Lambda^+}m_\lambda L(\lambda)\,,
		\]
		where $L(\lambda)$ denotes the finite-dimensional irreducible
		$\gf$-module of highest weight $\lambda$ and $m_\lambda$ is the multiplicity of $L(\lambda)$ in $M\otimes M$. Since $M$ is irreducible, we conclude from Schur's lemma that $\End_\gf(M)$ is a~division algebra. On the other hand, because $\gf$ is perfect (see~\cite[Corollary 5.2]{H})
		and $M$ is non-trivial, we have that $d:=\dim_\fbb M>1$, and therefore
		$\End_\fbb(M)\cong M_d (\fbb)$ has zero divisors. Consequently, we obtain that $\End_\gf(M)\subsetneqq\End_\fbb(M)$.

		Next, we compute $m_0$ by taking $\gf$-invariants on both sides, and as $M$ is self-dual, we then deduce the following inequality:
		\begin{eqnarray*}
			m_0 & = & \dim_\fbb(M\otimes M)^\gf=\dim_\fbb(M^*\otimes M)^\gf\\
			& = & \dim_\fbb\Hom_\fbb(M,M)^\gf=\dim_\fbb\End_\gf(M)\\
			& < & \dim_\fbb\End_\fbb(M)=\dim_\fbb(M\otimes M)\,.
		\end{eqnarray*}
		Hence, for dimension reasons, $M\otimes M$ must have a~non-trivial irreducible direct summand, and therefore we conclude from Lemma~\ref{criterion}\,(c)
		that $\K^{\bi} (\lf)$ is not a~Jordan ring. Finally, it follows from Lemma~\ref{jordanpowerasso}\,(a) that $\K^{\bi}(\lf)$ is also not alternative.
	\end{proof}

	\begin{remark}
		{Unfortunately, we do not know whether the Grothendieck ring in
		Theorem~\ref{nonalt} is power-associative or not. But a~straightforward computation using the Clebsch-Gordan formula shows that $\K^{\bi} (\slf_2
		(\cbb))$ is not power-associative, so we suspect this might be true more generally for every finite-dimensional non-zero semi-simple Leibniz algebra in characteristic zero. On the other hand, except for Proposition~\ref{nonassoK_0},
		Theorem~\ref{alt}, and Corollary~\ref{alt2}, we do not know anything about the associativity properties of $\K^{\bi}(\lf)$ for non-semi-simple
		Leibniz algebras $\lf$ in characteristic zero. Moreover, we know even less about these associativity properties in non-zero characteristics, but see Proposition~\ref{nonassoK_0}, Example~\ref{primchar}, and
		Proposition~\ref{nonalt2} below.}
	\end{remark}
	Be aware that, contrary to the situation in characteristic zero, the Killing form (or more generally, trace forms) of simple Lie algebras over a~field of non-zero characteristic can be zero (see Theorem 1.3 in Chapter 4 of~\cite{SF}). But nevertheless, similar to Theorem~\ref{nonalt}, one can deduce the following result which is valid for ground fields of arbitrary characteristic:
	\begin{proposition}\label{nonalt2}
		Let $\lf$ be a~non-zero Leibniz algebra whose canonical Lie algebra
		$\lf_{\lie}$ is finite dimensional simple and admits a~non-degenerate invariant bilinear form. Then the Grothendieck ring $\K^{\bi}(\lf)$ is neither alternative nor a~Jordan ring.
	\end{proposition}

	\begin{remark}
		{Note that the hypothesis of Proposition~\ref{nonalt2} is satisfied in many instances. Namely, for Lie algebras of classical type (see~\cite[Corollary 6.1]{B2} and~\cite[Theorems A and B]{G}), for graded Lie algebras of Cartan type (see Section 6 in Chapter~4 of~\cite{SF}), for (non-graded) Block algebras (see~\cite[Theorem 7]{B1}), and for Melikian algebras (see~\cite[Proposition~6.1]{PS}).}
	\end{remark}
	We conclude our paper by briefly discussing a~possible weak analogue of
	Theorem~\ref{communitalprod} for the Grothendieck ring $\K_{\weak}^{\bi}
	(\lf)$ of a~Leibniz algebra $\lf$.

	Recall that the \emph{pushout\/} of arbitrary commutative associative rings $R_1$ and $R_2$ with unities $1_{R_1}$ and $1_{R_2}$, respectively, is the tensor product $R_1\otimes_\zbb R_2$ with multiplication
	\[
		(r_1\otimes r_2 )(r_1 '\otimes r_2 '):=(r_1 r _1 ')\otimes(r_2 r_2 ')
	\]
	and unity $(1_{R_1},1_{R_2})$.

	Let $\lf$ be a~Leibniz algebra over a~field $\fbb$. Then we consider the following ring homomorphisms:
	\begin{gather*}
		\iota:\zbb\to\K(\lf_{\lie})\,,\,1\mapsto[F_0 ]\,, \\
		\varphi^s:\K(\lf_{\lie})\to\K_{\weak}^{\bi}(\lf)\,,\,[M]\mapsto[M^s ]\,, \\
		\varphi^a:\K(\lf_{\lie})\to\K_{\weak}^{\bi}(\lf)\,,\,[M]\mapsto[M^a ]\,.
	\end{gather*}
	Clearly, we have that $\varphi^s\circ\iota=\varphi^a\circ\iota$, and thus it follows from the universal property of the pushout that there is a~unique ring homomorphism
	\[
		\pi:\K(\lf_{\lie})\otimes_\zbb\K(\lf_{\lie})\to\K_{\weak}^{\bi}(\lf)
	\]
	such that $\varphi^s =\pi\circ\kappa^s$ and $\varphi^a
	=\pi\circ\kappa^s$, where $\kappa^s:\K(\lf_{\lie})\to\K(\lf_{\lie})
	\otimes_\zbb\K(\lf_{\lie})$ and $\kappa^a:\K(\lf_{\lie})\to\K(\lf_{\lie})
	\otimes_\zbb\K(\lf_{\lie})$ are the defining ring homomorphisms of the pushout. Note that one cannot even hope that $\pi$ is an epimorphism as it is not even an epimorphism of abelian groups
	(see our discussion of a~possible classification of finite-dimensional irreducible weak Leibniz bimodules in Section~\ref{weak}).

	On the other hand, $\Phi:\K^{\bi}(\lf)\to\K_{\weak}^{\bi}(\lf)$ defined by $[M^s ]\mapsto[M^s ]$ and $[M^a ]\mapsto[M^a ]$ on a~set of free generators of $\K^{\bi}(\lf)$ is a~monomorphism of abelian groups which clearly is not always compatible with the multiplications as the domain is not associative (at least in the case that $\dim_\fbb\lf_{\lie}
	<\infty$, see Proposition~\ref{nonassoK_0}), but the codomain is associative.

    
\end{document}